\documentclass[12pt]{amsart}
\usepackage{amsmath, amsthm, amsfonts, amssymb, color}
\usepackage{mathrsfs}
\usepackage{amstext,amsxtra}
\usepackage{txfonts} 
\usepackage{multirow}\usepackage{tabularx}
\usepackage{array}
\usepackage{caption}
\allowdisplaybreaks

\usepackage{pgf,tikz}
\usetikzlibrary{patterns}
\usetikzlibrary{arrows, arrows.meta}

\newlength{\hatchspread}
\newlength{\hatchthickness}
\newlength{\hatchshift}
\newcommand{\hatchcolor}{}
\tikzset{hatchspread/.code={\setlength{\hatchspread}{#1}},
         hatchthickness/.code={\setlength{\hatchthickness}{#1}},
         hatchshift/.code={\setlength{\hatchshift}{#1}},
         hatchcolor/.code={\renewcommand{\hatchcolor}{#1}}}
\tikzset{hatchspread=3pt,
         hatchthickness=0.4pt,
         hatchshift=0pt,
         hatchcolor=black}

\pgfdeclarepatternformonly[\hatchspread,\hatchthickness,\hatchshift,\hatchcolor]
   {custom north west lines}
   {\pgfqpoint{\dimexpr+2\hatchthickness}{\dimexpr+2\hatchthickness}}
   {\pgfqpoint{\dimexpr\hatchspread-2\hatchthickness}{\dimexpr\hatchspread-2\hatchthickness}}
   {\pgfqpoint{\dimexpr\hatchspread}{\dimexpr\hatchspread}}
   {
    \pgfsetlinewidth{\hatchthickness}
    \pgfpathmoveto{\pgfqpoint{0pt}{\dimexpr\hatchspread+\hatchshift}}
    \pgfpathlineto{\pgfqpoint{\dimexpr\hatchspread-0.15pt+\hatchshift}{0.15pt}}
    \ifdim \hatchshift < 0pt
      \pgfpathmoveto{\pgfqpoint{0pt}{\hatchshift}}
      \pgfpathlineto{\pgfqpoint{\dimexpr-0.15pt+\hatchshift}{0.15pt}}
    \fi
    \pgfsetstrokecolor{\hatchcolor}
    \pgfusepath{stroke}
   }

\newcommand{\RR}{\mathbb{R}}
\newcommand{\NN}{\mathbb{N}}
\newcommand{\Z}{\mathbb{Z}}
\newcommand{\CC}{\mathbb{C}}
\newcommand{\supp}{{\rm supp\,}}
\newcommand\LL{{\mathcal L}}
\renewcommand{\a}{{\alpha}}
\renewcommand{\b}{{\beta}}

\usepackage{color}

\newcommand{\red}[1]{\textcolor{red}{#1}}

\newtheorem{theorem}{Theorem}[section]

\newtheorem{proposition}{Proposition}[section]
\newtheorem{coro}{Corollary}[section]
\newtheorem{lemma}{Lemma}[section]

\newtheorem{remark}{Remark}[section]
\numberwithin{equation}{section}

\newcommand\R{\mathbb{R}}

\renewcommand{\r}{\right}
\renewcommand{\l}{\left}

\makeatletter
\def\l@subsection{\@tocline{2}{0pt}{2.5pc}{5pc}{}}
\makeatother

\begin{document}

\title{Bochner-Riesz means for Laguerre expansions}
\author{Himani Sharma, Adam Sikora, and  Sundaram Thangavelu}

\address{
	Himani Sharma, Formerly at the Department of Mathematics, Karlsruhe Institute of Technology (KIT), Karlsruhe-76131, Germany}
\email{himani.sharma@kit.edu}
\email{himani.sharma@alumni.anu.edu.au}
\address{
	Adam Sikora, Department of Mathematics, Macquarie University, NSW 2109, Australia}
\email{adam.sikora@mq.edu.au}
\address{
	Sundaram Thangavelu, National Science Chair, Department of Mathematics, Indian Institute of Science, Bangalore 560012, India\\
	and Honorary Professor, Department of Mathematics, University of Queensland, Brisbane, QLD 4072, Australia }
\email{veluma@iisc.ac.in}


\maketitle



\begin{abstract}
We investigate the $ L^p $ norm convergence of Bochner–Riesz means $ S_R^\delta f, \, f \in L^p(\R^+,d\mu_\b)$ associated to Laguerre expansions of type $\beta \ge -1/2$. For any given $ 1 < p < \infty$ and $ \beta > -1/2,$ we provide a complete characterisation of the order $\delta$ for which these means are uniformly bounded on $L^p$ spaces. A particularly noteworthy phenomenon arises for $-1/3 <\beta  <0$, where the admissible $(p,\delta)$  region differs substantially from those for $ -1/2 \leq \beta \leq -1/3$ or $\beta\ge0$. Interestingly, the convergence for $-1/3 < \beta <0$
is an interaction between the ranges $ -1/2 \leq \beta \leq -1/3$
and $\beta \ge 0$.
\end{abstract}

 \tableofcontents

\section{Introduction}
In this work, we investigate the Bochner-Riesz summability associated to  different kinds of Laguerre expansions. The study of such expansions has a long history, and one of the main aims in this article is to complete the picture by answering some questions that have remained open until now. More generally we study expansions in terms of a two-parameter family of Laguerre functions on suitable weighted  $ L^p $ spaces. These functions arise as eigenfunctions of certain second-order differential operators, namely  Laguerre differential operators, which generalise the simple harmonic oscillator Hamiltonian, also known as the Hermite operator.\\

In order to describe the problems we study, let us set up the notations first. For any $ \alpha > -1,$  Laguerre polynomials $L_{k}^{\alpha}$ of type $\alpha$ are defined by the equation  
\begin{equation*}
    e^{-x} x^{\alpha} L_{k}^{\alpha}(x)
    = \frac{1}{k!} \frac{d^{k}}{dx^{k}}
      \Big( e^{-x} x^{k+\alpha} \Big),\,\,\, k \in \mathbb N, x >0.
\end{equation*}
These polynomials arise  as eigenfunctions of the classical Laguerre differential operator, namely,
\begin{equation}\label{op-laguerre}
    x y''(x) + (\alpha + 1 - x) y'(x) + k\, y(x) = 0.
\end{equation}
The polynomials $\{L_{k}^{\alpha}\}_{k \geq 0}$ form an orthogonal basis in the weighted space
$
L^{2}\bigl(\R^+,\, e^{-x} x^{\alpha}\, dx \bigr),
$
for all  fixed $\alpha>-1$. As proved by Pollard \cite{Poll} the expansions in terms of $ L_k^\alpha$ are not amenable to Bochner-Riesz summability in the spaces $ L^{p}\bigl(\R^+,\, e^{-x} x^{\alpha}\, dx \bigr)$ for any $ p \neq 2.$  This has led to the study of expansions in terms of the Laguerre functions
\begin{equation}\label{NLagDef}
    \mathcal{L}_{k}^{\alpha}(x)
    = \Bigg( \frac{\Gamma(k+1)}{\Gamma(k+\alpha+1)} \Bigg)^{1/2}
      L_{k}^{\alpha}(x)\, e^{-x/2} x^{\alpha/2}.
\end{equation}
These functions form a complete orthonormal system in $L^2(\R^+,dx)$ and a natural problem is to study partial sums and Bochner-Riesz means associated to expansions in terms of  $  \mathcal{L}_{k}^{\alpha}.$ These expansions are known as standard Laguerre expansions in the literature. \\

In 1965, motivated by the results of Pollard and Wing, Askey and Wainger undertook the study of Hermite and Laguerre expansions in  \cite{AW}. Under the assumption that $ \alpha \geq 0, $ they considered the partial sums 
$$ \sigma_N f(x) = \sum_{k=0}^N ( f,   \mathcal{L}_{k}^{\alpha})\,   \mathcal{L}_{k}^{\alpha}(x),\,\,\, f \in L^p(\R^+, dx).$$
By using very precise asymptotic behaviour of Laguerre functions and weighted norm inequalities for the Hilbert transform and Hardy's averaging operators, they  managed to prove a somewhat surprising result: the partial sums $ \sigma_N $ are uniformly bounded on $ L^p(\R^+,dx) $ if and only if $ 4/3 < p < 4.$ (They also proved a similar result for expansions in terms of Hermite functions.)  Their work demonstrated that appropriate renormalisation and oscillatory estimates can restore $L^{p}$ convergence in situations where the weighted theory fails, and their methods have since influenced the broader development of multiplier theorems on orthogonal polynomial systems.\\

The negative result of Pollard and the positive result of Askey and Wainger triggered quite a lot of interest in the study of Bochner-Riesz  and Ces{\`a}ro means associated to the Hermite and Laguerre expansions. For standard Laguerre expansions, fundamental contributions were made by Muckenhoupt \cite{Muck1, Muck2}, who established mean convergence theorems for both Laguerre and Hermite expansions across a broad range of parameters. His analysis introduced techniques inspired by what later became the theory of $A_{p}$ weights, providing precise sufficient conditions for the weighted $L^{p}$-boundedness of partial sums and related operators. These results continue to form the analytic foundation for modern studies of spectral multipliers and convergence questions for Laguerre expansions. \\

The  Ces{\`a}ro means of order $ \delta \geq 0,$ associated with standard Laguerre expansions are defined by
$$  \sigma_N^\delta f(x) = \frac{1}{A_N^\delta} \sum_{k=0}^N  A_{N-k}^\delta( f,   \mathcal{L}_{k}^{\alpha})\,   \mathcal{L}_{k}^{\alpha}(x),\,\,\, f \in L^p(\R^+, dx)$$
where  $ A_N^\delta $ are  the generalised binomial coefficients  defined by the equation
$$
A_N^\delta = \frac{\Gamma(N+\delta+1)}{\Gamma(N+1)\Gamma(\delta+1)}.
$$
The uniform boundedness of these means on $ L^p(\R^+,dx) $ was studied by Poiani \cite{Poi72} when $ \delta = 1$  and by Markett \cite{mark1} in the general case. 
Under the assumption that $ \alpha \geq 0,$ Markett proved  in a series of papers \cite{mark1, mark2, mark3} the most definitive results on the $ L^p$ boundedness of $ \sigma_N^\delta.$ In particular, he proved that the critical index for $ L^1 $ summability is $ \delta= 1/2.$ \\

Markett proved his positive results by making use of  a convolution structure called the Laguerre convolution. In this context, he was led to consider expansions in terms of the functions 
\begin{equation}\label{LagFunc}
    \varphi_k^\alpha(x) := \Bigg( \frac{2\Gamma(k+1)}{\Gamma(k+\alpha+1)} \Bigg)^{1/2} L_k^\alpha(x^2) e^{-x^2/2} .
\end{equation}
These functions $\{\varphi_k^\a\}_{k \ge 0}$ form an orthonormal basis for the Hilbert space $L^2(\R^+, x^{2\alpha+1} dx)$. Markett studied the Ces{\`a}ro means associated to expansions in terms of  these functions. For $ f \in L^2(\R^+, x^{2\alpha+1} dx)$ these means are defined by
$$  \sigma_N^\delta f(x) = \frac{1}{A_N^\delta} \sum_{k=0}^N  A_{N-k}^\delta \langle f,   \varphi_{k}^{\alpha}\rangle\,   \varphi_{k}^{\alpha}(x).$$
where we have defined
$$\langle f,   \varphi_{k}^{\alpha}\rangle\, = \int_0^\infty f(x)\, \varphi_k^\alpha(x)\, x^{2\alpha+1}\, dx.$$
Under the assumption that $ \alpha \geq 0,$ he proved that   $ \delta = \alpha+1/2 $ is the critical index for the $ L^1 $ summability and for $ 0 < \delta \leq \alpha+1/2,$  the above Ces{\`a}ro means are uniformly bounded on $ L^p(\R^+, x^{2\alpha+1}\, dx) $ if and only if $ \frac{4\alpha+4}{2\alpha+3+2\delta} < p <  \frac{4\alpha+4}{2\alpha+1-2\delta}.$ \\

Markett noticed that the assumption $ \alpha \geq 0 $ cannot be dispensed with since the proof depends on the Young's inequality for the Laguerre convolution which is known only for $ \alpha \geq 0.$ It can be shown that the boundedness of the Laguerre translation on $ L^1(\R^+, x^{2\alpha+1}dx ) $ cannot be true for $ -1/2 \leq \alpha \leq -1/3 $ but for $ \alpha >-1/3 $ it is still an open problem. Understanding the behaviour of Ces{\`a}ro/Bochner-Riesz means  associated to expansions in terms of $ \varphi_k^\alpha, $ for all values of $ \alpha \geq -1/2 $ was our main concern when we started this work. \\

Expansions in terms of $ \varphi_k^\alpha $ are  known as Laguerre expansions of convolution type in the literature. These expansions are also intimately connected with Hermite and special Hermite expansions.
The one dimensional Hermite functions $ h_k(x) $ are defined by the equation
$$ h_k(x) = \Big( \sqrt{\pi} 2^k k! \Big)^{-1/2} (-1)^k \, e^{\frac{1}{2} x^2}\, \frac{d^k}{dx^k}e^{-x^2} $$
and they form an orthonormal basis for $ L^2(\R).$ In higher dimensions, the Hermite functions $ \Phi_\mu(x), \mu \in \mathbb N^n $ are defined by taking tensor products of the one dimensional functions. These functions are eigenfunctions of the Hermite operator $ H=-\Delta+|x|^2 $ with eigenvalues $ (2|\mu|+n).$ The connection between Hermite expansions and Laguerre expansions of type $(n/2-1)$ is given by the following relation. For any radial function $ f(x) = g(|x|) $ in $ L^2(\R^n) $ we have
$$  \sum_{|\mu|=2k} (f, \Phi_\mu) \Phi_\mu(x) =c_n\, \langle g, \varphi_k^{n/2-1} \rangle \, \varphi_k^{n/2-1}(|x|) $$ 
and $ (f, \Phi_\mu) =0 $ for odd values of $ |\mu|,$  see e.g. \cite{Th2} or \cite[Section 3.4 ]{ThBook}. Similarly expansions in terms of $ \varphi_k^{n-1} $ are connected to special Hermite expansions, i.e. expansions in terms of eigenfunctions of the special Hermite operator
$$ L = -\Delta_{\mathbb C^n}+\frac{1}{4}|z|^2-i \sum_{j=1}^n \big( x_j \frac{\partial}{\partial y_j}- y_j \frac{\partial}{\partial x_j}\big)$$ 
on $\mathbb C^n ,$ also known as the twisted Laplacian. For more details of this connection  see \cite[Section 6.1]{ThBook} where a motivation for the Laguerre convolution can also be found.\\

Even at an elementary level, the one dimensional Hermite functions $ h_k $ are related to $ \varphi_k^\alpha $ for $ \alpha= -1/2, 1/2.$ Indeed,
$$ h_{2k}(x) = (-1)^k  \varphi_k^{-1/2}(x),\,\,\, h_{2k+1}(x) = (-1)^k  \varphi_k^{1/2}(x)\, x .$$
These relations motivated Markett to consider expansions in terms of the functions
$$ h_k^\alpha(x) = \mathcal{L}_k^\alpha(x^2) (2x)^{1/2} = \sqrt{2}\,\varphi_k^\alpha(x)\, x^{\alpha+1/2}.$$
These functions form an orthonormal basis for $ L^2(\R^+) $ and expansions in terms of them are known as Laguerre expansions of Hermite type. By studying the behaviour of the  $ L^1(\R^+) $ norms of the partial sum operators, Markett observed that Laguerre expansions of Hermite type has a better behaviour than the standard Laguerre expansions. This investigation led him to conjecture that the critical index for the $ L^1 $ summability for the Hermite type expansions satisfies $ 1/6 \leq \delta \leq 1/2.$ This conjecture turned out to be true as shown in \cite{Th90} for $\a\ge1/2$, and for other cases in \cite[Section 6.3]{ThBook}.\\ 

Returning to Laguerre expansions of convolution type, namely in terms of the functions $ \varphi_k^\alpha$,  let us denote  the critical index for  $ L^1 $ summability by $ \delta(\alpha).$  Since Laguerre expansions of type $ \alpha = -1/2 $ are subsumed by Hermite expansions, from the results  proved in \cite{Th89}  we infer that $ \delta(-1/2) = 1/6 $  whereas from the works of Markett we also have $ \delta(\alpha) =\alpha+1/2 $ for $ \alpha \geq 0.$   Based on very precise estimates  on the $ L^p $  norms of the Laguerre functions $ \varphi_k^\alpha $ with respect to the measure $ d\mu_\alpha = x^{2\alpha+1} dx$ we are led to the conjecture that  the critical index satisfies $ \delta(\alpha) = 1/6,$ for $-1/2 \leq \alpha \leq -1/3,$ and $ \delta(\alpha) = \alpha+1/2$ for $ \alpha \geq -1/3.$  We prove this conjecture and more in this article.\\

Since the Laguerre polynomials $ L_k^\alpha(x) $ satisfy the eigenfunction equation \eqref{op-laguerre}, it follows that the standard Laguerre functions are eigenfunctions of the self-adjoint operator
$$
\mathcal{L}_{\alpha}= -x \frac{d^{2}}{dx^{2}}- (\alpha + 1 - x) \frac{d}{dx}+ \frac{x}{4} - \frac{\alpha^{2}}{4x},
$$
with eigenvalues $\lambda_{k} = (2k + \alpha+1)/2$. By standard computations, one can easily check  that $\varphi_k^\a$ are eigenfunctions of the differential operator defined by 
\begin{equation}\label{RadialHO}
    \Delta_\alpha = -\frac{d^2}{dx^2} - \frac{2\alpha + 1}{x} \frac{d}{dx} + x^2,
\end{equation}
with eigenvalues $4k+2\a+2$. That is,
\begin{equation*}
    \Delta_\alpha \varphi_k^\alpha(x) = (4k + 2\alpha + 2)\, \varphi_k^\alpha(x).
\end{equation*}
We observe that for $n = 2\alpha + 2$, the operator $\Delta_\alpha$ corresponds to the radial part of the Hermite operator  $-\Delta + |x|^2$ acting on $L^2(\mathbb{R}^n)$.  For this reason, $ \Delta_\alpha $ may be termed as fractional  harmonic oscillator. The isometry
\begin{equation*}
    M_{1/2+\alpha} f(x) = x^{1/2 + \alpha} f(x)
\end{equation*}
provides a unitary equivalence between $L^2(\R^+, dx)$ and $L^2(\R^+, x^{2\alpha + 1} dx)$. Under this transformation, $\Delta_\alpha$ becomes unitarily equivalent to the Schrödinger-type operator
\begin{equation}\label{RadialHO2}
    \widetilde{\Delta}_\alpha = -\frac{d^2}{dx^2} + \frac{\alpha^2 - 1/4}{x^2} + x^2,
\end{equation}
acting on $L^2(\R^+, dx)$, i.e.,
\begin{equation*}
    M_{1/2+\alpha}^{-1} \widetilde{\Delta}_\alpha M_{1/2+\alpha} f(x) = \Delta_\alpha f(x).
\end{equation*}
The eigenfunctions of the operator $\widetilde{\Delta}_\alpha$  are provided by Laguerre  functions of Hermite type, namely 
\begin{equation}\label{eqDelta2}
    \widetilde{\Delta}_\alpha h_k^\alpha(x) = (4k + 2\alpha + 2)\, h_k^\alpha(x).\\
\end{equation}

An interesting feature of these operators is their behaviour under the transformation $\alpha\mapsto-\a$. For the operator $\Delta_\a$, both the differential expression and the underlying measure $x^{2\a+1}dx$ depend explicitly on $\a$, so replacing $\a$ by $-\a$ leads to a genuinely different operator acting on a different space, and consequently different spectral properties.

In contrast, the Schrödinger-type operator $\tilde{\Delta}_\a$ is formally invariant under the change $\a\mapsto-\a$, since the potential depends only on $\a^2$. However, this symmetry is broken at the level of boundary conditions. A standard analysis shows that near $x=0$, solutions to $\tilde{\Delta}_\a f=\lambda f$ behave like
$$ f(x)\sim x^{\a+1/2} \quad \text{ or } \quad x^{-\a+1/2}.$$
For $|\a|\ge 1$, only one of the above two behaviours is square-integrable, and thus the self-adjoint realisation selects the boundary condition corresponding to either $f(x)=O(x^{\a+1/2})$ or $f(x)=O(x^{-\a+1/2})$ as $x\to 0$. Under the transformation $\a\mapsto-\a$, these two asymptotic behaviours are interchanged, and hence the admissible boundary condition at the origin changes accordingly.

Thus, although the formal differential expression of $\tilde{\Delta}_\a$ remains unchanged, the associated self-adjoint operator, and therefore its spectral resolution, is affected through the choice of boundary condition at the singular point $x=0$. Throughout the paper, we always consider these natural self-adjoint realisations of the operators.\\

The form of the operator $  \widetilde{\Delta}_\alpha $ suggests that we consider the following two parameter family of operators
defined by
\begin{equation}\label{MainOp}
    \Delta_{\alpha,\beta} := -\frac{d^2}{dx^2} - \frac{2\beta + 1}{x} \frac{d}{dx} + \frac{(\alpha - \beta)(\alpha + \beta)}{x^2} + x^2,
\end{equation}
acting on the spaces $L^p(\R^+, x^{2\beta + 1} \, dx)$, where $\alpha, \beta \geq -\tfrac{1}{2}$. The functions $ \psi_k^{\alpha,\beta} $ defined by
\begin{equation}\label{ab}
    \psi_k^{\alpha,\beta}(x) := \varphi_k^{\alpha}(x) \, x^{\alpha - \beta} = \mathcal{L}_k^{\alpha}(x^2) \, x^{-\beta},
\end{equation}
turn out to be eigenfunctions of the operator $ \Delta_{\alpha,\beta}.$ More precisely,
\begin{equation*}
     \Delta_{\alpha,\beta}\, \psi_k^{\alpha,\beta}(x)=(4k+2\a+2)\psi_k^{\alpha,\beta}(x).
\end{equation*}
It is clear that the family $ \psi_k^{\alpha,\beta}, k \geq 0 $ forms an orthonormal basis for  $L^2(\R^+, x^{2\beta + 1} \, dx).$  We are interested in Bochner-Riesz means and spectral multipliers associated to expansions in terms of $ \psi_k^{\alpha,\beta}.$ Note that this covers all three kinds of Laguerre expansions since $ \psi_k^{\alpha,0}(x) = \mathcal{L}_k^\alpha(x^2), \psi_k^{\alpha,\alpha}(x) = \varphi_k^\alpha(x) $ and  $ \psi_k^{\alpha,-1/2}(x) = h_k^\alpha(x).$\\

It is interesting to note that  $\Delta_{\pm\alpha} = \Delta_{\alpha,\pm\alpha} $ and $ \widetilde{\Delta}_\alpha = \Delta_{\alpha, -1/2}.$
Moreover, the unitary map
\begin{equation}\label{unitary1}
    M_{1/2 + \beta} \, \Delta_{\alpha,\beta} \, M^{-1}_{1/2 + \beta} f(x) = \widetilde{\Delta}_\alpha f(x)
\end{equation}
establishes equivalence between $\Delta_{\alpha,\beta}$ and $\widetilde{\Delta}_\alpha$. Additionally, we have the identity
\begin{equation}\label{unitary2}
    M^{-1}_{\alpha - \beta} \, \Delta_{\alpha,\beta} \, M_{\alpha - \beta} f(x) 
    = M^{-1}_{1/2 + \alpha} \, \widetilde{\Delta}_\alpha \, M_{1/2 + \alpha} f(x) 
    = \Delta_\alpha f(x),
\end{equation}
illustrating that all three operators are unitarily equivalent.\\


For any fixed $ \alpha \geq -1/2,$ all the operators $ \Delta_{\alpha,\beta} $ have the same spectrum, namely  $ \lambda_k = \lambda_k(\alpha) =(4k+2\a+2), \, k \geq 0.$  
The Bochner-Riesz means of  order $ \delta \geq 0 $  associated to expansions in terms of $ \psi_k^{\a,\b} $  are defined by 
\begin{equation}\label{BRmeans}  S^\delta_R(\Delta_{\alpha,\beta})f(x) = \sum_{\lambda_k \leq R} (1-\lambda_k/R)_+^\delta \langle f, \psi_k^{\alpha,\beta} \rangle_\beta\, \psi_k^{\alpha,\beta}(x).
\end{equation}
In the above,  $ \langle f, \psi_k^{\alpha,\beta}\rangle_\beta$ are the Laguerre coefficients  of the function $ f $  in $ L^2(\R^+, x^{2\beta+1} dx)$ defined by
$$ \langle f, \psi_k^{\alpha,\beta} \rangle_\beta = \int_0^\infty  f(x)\,  \psi_k^{\alpha,\beta}(x)\, x^{2\beta+1}\, dx.$$
Sometimes we simply write $ S_R^\delta $ when the operator under consideration is fixed.  It is known that the uniform boundedness of $ S_R^\delta ,\, R >0$ is equivalent to that of the Ces{\`a}ro means $ \sigma_N^\delta, N \in \mathbb N.$
For any given $ \alpha, \beta \geq -1/2,$ and $ 1 \leq p \leq \infty,$ we are interested in finding necessary and sufficient conditions on $ \delta $ so that the operators $S_R^\delta(\Delta_{\alpha,\beta}) $ are uniformly bounded on $ L^p(\R^+, x^{2\beta+1} dx) $ for a given $p.$\\

The Bochner-Riesz means are a special case of spectral multipliers. Given a function  $ F $ on the spectrum,  the spectral multiplier $ F(\Delta_{\alpha,\beta}) $ is defined on $ L^2(\R^+, x^{2\beta+1}dx) $ by 
$$ F(\Delta_{\alpha,\beta})f(x) = \sum_{k=0}^\infty F(\lambda_k)\, \langle f, \psi_k^{\alpha,\beta} \rangle_\beta\, \psi_k^{\alpha,\beta}(x). $$ 
Note that $ S_R^\delta(\Delta_{\a,\b})= s_R^\delta(\sqrt{\Delta_{\a,\b}}) $ if we take  $ s_R^\delta(\lambda) = (1-\lambda^2/R)_+^\delta.$  As we will see later, this relation is crucial in handling certain parts of $ S_R^\delta$  where the finite propagation speed of the wave equation associated to $ \Delta_{\a,\b} $ will be used.\\

Both Askey–Wainger and Muckenhoupt relied on detailed asymptotic estimates for Hermite and Laguerre functions across different spatial regimes, together with weighted norm inequalities for the Hilbert transform and the Hardy operator. The weights required in their arguments were elementary power weights, so the general $A_p$ theory was not needed. Indeed, the
$A_p$ theory itself was developed by Muckenhoupt only later, beginning with his 1972 paper and subsequently in joint work with Hunt and Wheeden on the Hilbert transform on the circle.

\section{Preliminaries}
In this section we collect some relevant estimates, both pointwise and in norm, on  various  Laguerre functions and their differences. These estimates are well-known and   provide key ingredients for the proofs of our main theorems. We begin with pointwise estimates.\\

We recall the following estimates on the Laguerre functions $\LL^\a_k$ from \cite{Muck1}, see also \cite{AW}, which show that for $\a>-1$, there are positive constants, $C$ and $\gamma$, independent of $x$ and $k$ such that for all integers $k\ge0$,
\begin{equation}\label{LaguerreBounds}
    |\LL^\a_{k}(x)|\leq C\begin{cases}
        x^{\a/2}\lambda_k^{\a/2},  \quad &0\le x\le 1/\lambda_k,\\
      x^{-1/4}\lambda_k^{-1/4}, \quad &1/\lambda_k< x\le \lambda_k/2 ,\\
     \lambda_k^{-1/4}(\lambda_k^{1/3}+|x-\lambda_k|)^{-1/4}, \quad &\lambda_k/2< x\le 3\lambda_k/2, \\
      e^{-\gamma x}, \quad &3\lambda_k/2<x,
    \end{cases}
\end{equation}
where $\lambda_k=4k+2\a+2.$ Similarly for $k\ge0$
\begin{equation}\label{LagDiffBounds}
    |\LL^\a_{k}(x)-\LL_{k+2}^\a(x)|\leq C\begin{cases}
        x^{\a/2}\lambda_k^{\a/2-1}, \quad &0\le x\le 1/\lambda_k,\\
        x^{1/4}\lambda_k^{-3/4}, \quad &1/\lambda_k< x\le \lambda_k/2 ,\\
        \lambda_k^{-3/4}(\lambda_k^{1/3}+|x-\lambda_k|)^{1/4}, \quad &\lambda_k/2< x\le 3\lambda_k/2, \\
        e^{-\gamma x}, \quad &3\lambda_k/2<x.
    \end{cases}
\end{equation}
Note that if $k=0$, and $\a$ is negative, then $\lambda_k$ should be taken equal to 2 in \eqref{LaguerreBounds}.\\

Next, 
from the above estimates we get the following upper bounds for $\psi^{\a,\b}_k(x)$
\begin{equation}\label{psi}
    |\psi^{\a,\b}_k(x)|=|\LL^\a_{k}(x^2)x^{-\b}|\leq C \mathcal{G}^{\a,\b}_k(x),
\end{equation}
where 
\begin{equation*}
    \mathcal{G}^{\a,\beta}_k(x)=\begin{cases}
        x^{\a-\b}\lambda_k^{\a/2},  \quad &0\le x^2\le 1/\lambda_k,\\
      x^{-\b-1/2}\lambda_k^{-1/4}, \quad &1/\lambda_k< x^2\le \lambda_k/2 ,\\
     x^{-\b}\lambda_k^{-1/4}(\lambda_k^{1/3}+|x^2-\lambda_k|)^{-1/4}, \quad &\lambda_k/2< x^2\le 3\lambda_k/2, \\
    e^{-\gamma x}, \quad &3\lambda_k/2<x^2
    \end{cases}
\end{equation*}
with $\lambda_k=4k+2\a+2.$ Similarly for $k\ge0$
\begin{equation}
   |\psi^{\a,\b}_k(x)-\psi^{\a,\b}_{k+2}(x)|= |(\LL^\a_{k}(x^2)-\LL_{k+2}^\a(x^2))x^{-\b}|\leq C\mathcal{H}^{\a,\b}_k(x),
\end{equation}
where
\begin{equation*}
  \mathcal{H}^{\a,\b}_k(x)=\begin{cases}
        x^{\a-\b}\lambda_k^{\a/2-1}, \quad &0\le x^2\le 1/\lambda_k,\\
        x^{-\b+1/2}\lambda_k^{-3/4}, \quad &1/\lambda_k< x^2\le \lambda_k/2 ,\\
        x^{-\b}\lambda_k^{-3/4}(\lambda_k^{1/3}+|x^2-\lambda_k|)^{1/4}, \quad &\lambda_k/2< x^2\le 3\lambda_k/2, \\
        e^{-\gamma x}, \quad &3\lambda_k/2<x^2.
    \end{cases}
\end{equation*}

Having collected all the relevant pointwise estimates, we now turn our attention towards $L^p $ norms of various Laguerre functions. All the estimates are deduced from the basic lemma proved by Markett \cite{mark2}, see
also  \cite[Lemma 1.5.4]{ThBook} which gives weighted estimates on the $L^p$ norms of $\LL^{\a+\beta}_k(x).$ \\

\begin{lemma}\label{Lag-p-norm}
Let $\alpha +\beta > -1$ and $\alpha >-2/p$, $1 \le p \le \infty$. Then, as $k \to \infty$ the following estimates are valid:
\begin{align*}
(i) &&\|{\mathcal L}_k^{\alpha+\beta}(x) x^{-\beta/2} \|_p \sim k^{1/p-1/2-\beta/2}
\quad \mbox{if} \quad \beta \le 2/p-1/2
\\
(ii) &&\|{\mathcal L}_k^{\alpha+\beta}(x) x^{-\beta/2} \|_p \sim k^{\beta/2-1/p}
\quad \mbox{if} \quad \beta > 2/p-1/2
\end{align*}
for p in the interval $1 \le p \le 4$. For $p>4$
\begin{align*}
(iii) &&\|{\mathcal L}_k^{\alpha+\beta}(x) x^{-\beta/2} \|_p \sim k^{\frac{-1}{3}+ \frac{1}{3p}-\beta/2  }
\quad \mbox{if} \quad \beta \le  \frac{4}{3p}-1/3
\\
(iv) &&\|{\mathcal L}_k^{\alpha+\beta}(x) x^{-\beta/2} \|_p \sim k^{\beta/2-1/p}
\quad \mbox{if} \quad \beta > \frac{4}{3p}-1/3.
\end{align*}
	\end{lemma}
\begin{proof}
See \cite[Lemma 1.5.4]{ThBook}; see also \cite{Muck1} for related estimates.
\end{proof}

As a consequence of this lemma,  it is not difficult to prove very precise estimates for  the $L^p$ norms of the eigenfunctions $\psi^{\a,\b}_k$ for different values of $\a$ and $\b$. For the convenience of the readers we first consider the case $ \alpha= \beta.$ Recall that when $\a=\b$, $\psi^{\a,\b}_k=\varphi^\a_k$. \\
  
  \begin{proposition}\label{basiclemma} For $ \alpha > -1,$  let $ \beta^\prime = 2\alpha(1/2-1/p).$ Then we have the following estimates: for $ 1 \leq p \leq 4$ we have
  $$ \| \varphi_k^{\alpha} \|_{L^p(d\mu_\alpha)} = O( k^{1/p-1/2-\beta^\prime/2}),\,\,\, \beta^\prime < 2/p-1/2,$$
   $$ \| \varphi_k^\alpha \|_{L^p(d\mu_\alpha)} = O( k^{-1/p+\beta^\prime/2}),\,\,\, \beta^\prime > 2/p-1/2.$$
   And for $ 4 < p < \infty $ we have
    $$ \| \varphi_k^\alpha \|_{L^p(d\mu_\alpha)} = O( k^{-1/3+1/(3p)-\beta^\prime/2}),\,\,\, \beta^\prime \leq 4/(3p)-1/3 ,$$
     $$ \| \varphi_k^\alpha \|_{L^p(d\mu_\alpha)} = O( k^{-1/p+\beta^\prime/2}),\,\,\, \beta^\prime > 4/(3p)-1/3.$$
  \end{proposition}
\begin{proof} 
  In view of the relation between $ \varphi_k^\alpha $ and  $ \mathcal{L}_k^\alpha$ we have the equation
  $$ \int_0^\infty |\varphi_k^\alpha(x)|^p d\mu_\alpha(x) = \int_0^\infty |\mathcal{L}_k^\alpha(x) x^{-\alpha /2} |^p\, x^\alpha dx.$$
  Therefore, we can rewrite the norm of $ \varphi_k^\alpha $ in $ L^p(\R^+, d\mu_\alpha) $ as
  $$  \int_0^\infty |\varphi_k^\alpha(x)|^p d\mu_\alpha(x)  = \int_0^\infty |\mathcal{L}_k^{\alpha^\prime+\beta^\prime}(x) x^{-\beta^\prime/2} |^p\, dx,$$
  where $ \alpha^\prime = 2\alpha/p$ and $ \beta^\prime = \alpha - \alpha^\prime = \alpha - 2\alpha/p = \alpha(1-2/p).$ Thus we can directly appeal to Lemma \ref{Lag-p-norm}. Note that when $ \alpha > -1,$ the condition $ \alpha^\prime > -2/p $ is satisfied.
  \end{proof} 

The estimates stated in the above proposition are in terms of $ \beta^\prime $ which is a function of $ \alpha.$ It is convenient to restate the estimates in terms of $ \alpha $ for the  ease of using them when needed. In what follows we only consider the case $ -1/2 \leq \a \leq 0$  to fix the ideas. The general case will be treated in a later proposition.\\

Consider $ \beta^\prime- (2/p-1/2) = \alpha( 1-2/p)-(2/p-1/2) = (2\alpha+1)/2 -2(\alpha+1)/p. $ Since $ (2\alpha+1) < (\alpha+1) $ for $ -1/2 \leq \alpha <0 $ we note that for $ 1 \leq p \leq 4, $ we have $ 1 \leq p < 4 (\alpha+1)/(2\alpha+1) $ so that  $ \beta^\prime < (2/p-1/2).$ This leads to the estimate
 $$ \| \varphi_k^\alpha \|_{L^p(d\mu_\alpha)} = O( k^{1/p-1/2-\beta^\prime/2}) = O(k^{(\alpha+1)(1/p-1/2)}),\,\,\,  1 \leq p \leq 4. $$
 When $  \alpha> -1/3,$ we note that $ (3\alpha+2) > 2(3\alpha+1) $ and hence we can consider $ 4 \leq p < 2(3\alpha+2)/(3\alpha+1).$  Consider
 $$ \beta^\prime - (4/(3p)-1/3) = (\alpha+1/3)-2(\alpha+2/3)/p = (1/3) \left((3\alpha+1)-2(3\alpha+2)/p\right) $$
 which is negative for  $ 4 \leq p < 2(3\alpha+2)/(3\alpha+1).$ In this case, we have the estimate
  $$ \| \varphi_k^\alpha \|_{L^p(d\mu_\alpha)} = O( k^{-1/3+1/(3p)-\beta^\prime/2}) = O(k^{(\alpha+1/3)/p-\alpha/2-1/3)})= O( k^{\alpha/2-(\alpha+1/3)/p^\prime} ).$$
  Still assuming $ \alpha > -1/3 $ but for $ p > 2(3\alpha+2)/(3\alpha+1) $ we have the estimate
  $$ \| \varphi_k^\alpha \|_{L^p(d\mu_\alpha)} = O( k^{\alpha/2-(\alpha+1)/p}) = O(k^{(\alpha+1)/p^\prime-\alpha/2-1}).$$
  It remains to consider the case $ -1/2 \leq \alpha \leq -1/3, p >4.$ In this case $ \beta^\prime - (4/(3p)-1/3) $ is always negative, and hence we have
   $$ \| \varphi_k^\alpha \|_{L^p(d\mu_\alpha)} = O( k^{-1/3+1/(3p)-\beta^\prime/2}) = O(k^{(\alpha+1/3)/p-\alpha/2-1/3})= O( k^{\alpha/2-(\alpha+1/3)/p^\prime} ).$$

As a corollary of the above estimates, we can prove the following result which gives a necessary condition for the $ L^p $ summability of the  Ces{\`a}ro means.\\
   
   \begin{coro}\label{basicprop} Assume that $ -1/2 \leq \alpha \leq -1/3.$ For $ 1 \leq p \leq 4/3 $ we have
   $$  \| \varphi_k^\alpha \|_{L^p(d\mu_\alpha)} \,  \| \varphi_k^\alpha \|_{L^{p^\prime}(d\mu_\alpha)} = O(k^{-1/2+2/(3p)}).$$
   Assume that  $ 0 > \alpha > -1/3 $ so that $ 1 <(6\alpha+4)/(3\alpha+3) < 4/3.$  For $ 1 \leq p < (6\alpha+4)/(3\alpha+3)$ 
    $$  \| \varphi_k^\alpha \|_{L^p(d\mu_\alpha)} \,  \| \varphi_k^\alpha \|_{L^{p^\prime}(d\mu_\alpha)}  = O(k^{2(\alpha+1)(1/p-1/2)-1/2}).       $$
    For $ (6\alpha+4)/(3\alpha+3) \leq p < 4/3$ we have
    $$  \| \varphi_k^\alpha \|_{L^p(d\mu_\alpha)} \,  \| \varphi_k^\alpha \|_{L^{p^\prime}(d\mu_\alpha)}  = O(k^{-1/2+2/(3p)}).       $$
    \end{coro}
    
    \vskip0.15truein
  
   For the uniform boundedness of the  Ces{\`a}ro means  $ \sigma_N^\delta $ on $ L^p(d\mu_\alpha)$ we get the following necessary condition on $ \delta.$ \\
   (i) When $ -1/2 \leq \alpha \leq -1/3, 1 \leq p <4/3 $  we need $ \delta  \geq -1/2+2/(3p).$ 
   (ii) When $ \alpha >-1/3 $  we have two cases to consider: If $1 \leq p < (6\alpha+4)/(3\alpha+3),$ then we need $ \delta \geq 2(\alpha+1)(1/p-1/2)-1/2.$ If $ (6\alpha+4)/(3\alpha+3) \leq p < 4/3,$ then we need $ \delta \geq -1/2+2/(3p).$\\
   
   Note that for  $ \alpha >-1/3 $ we have  two necessary conditions depending on $ p $ which coincide at $p = (6\alpha+4)/(3\alpha+3).$\\

If we denote by  $P^\a_k$ the projection operator corresponding to $\Delta_\a$ on $L^p(\R^+,x^{2\a+1}dx)$, that is 
\begin{equation}
   P^{\a}_kf:=\langle f,\varphi_k^{\a} \rangle_\a \varphi_k^{\a},
\end{equation}
where $\langle f,g \rangle_\a=\int_{\R^+}f(x) \bar{g}(x)x^{2\a+1}dx$ then the  figures below give the $L^p-L^2$ norm estimates of $P^\a_k$ for different ranges of $\a$. These figures are plotted on $(\a,\frac1{p'})$ plane, where $\frac1p+\frac1{p'}=1$. By duality, the $L^p\to L^2$ bounds for the projections $P^{\a}_k$ obtained here are equivalent to $L^2\to L^{p'}$ estimates. Therefore, in this form, they can be directly compared with the $L^2\to L^p$ spectral projection estimates for the Hermite operator on $\R^n$ established by Koch and Tataru \cite{KT} when $\a\ge 0$.

\begin{figure*}[!h]
    \centering
        \begin{tikzpicture}[scale=8]
\draw (0,0.55) -- (0,-0.25)-- (0,0) -- (1.10,0) --(1,0) -- (1,0.55);
\draw (0.25,0.01)--(0.25,-0.01); \draw (0.75,0.01)--(0.75,-0.01);
\draw (0.5,0.01)--(0.5,-0.01);
 \draw (3/8,0.01)--(3/8,-0.01);  \draw (5/8,0.01)--(5/8,-0.01);
\draw (0.99,1/2)--(1.01,1/2);\draw (-0.01,0.25)--(0.01,0.25);\draw (-0.01,1/2)--(0.01,1/2);
\draw (3/8,-0.125)--(1,1/2); \draw (3/8,-0.125)--(0,0.25);
\draw (-0.01,-0.125)--(0.01,-0.125);
\draw (0,0) node[left] {$0$} (1/4,0) node[below] {$\frac{2\a}{4\a+4}$} (3/8,0) node[above] {$\frac{2\a+1}{4\a+4}$} (1/2,0) node[below] {$\frac12$}
(5/8,0) node[below] {$\frac{2\a+3}{4\a+4}$} (3/4,0) node[below] {$\frac{3}{4}$} (1,0) node[below] {$1$} (1.1,0) node[below] {$\frac{1}{p'}$};
\draw (1,1/2) node[right] {$\frac{\alpha+1}{2}$} (0,-0.125) node[left] {$\frac{-1}{4}$} (0,0.25) node[left] {$\frac{\alpha}{2}$}; 
\end{tikzpicture}
    \caption{$\a \ge 0 $ }
    \label{FigLem2.2}
\end{figure*}

\begin{figure}[!h]
    \centering
    \begin{minipage}[b]{0.45\textwidth}
    \begin{tikzpicture}[scale=6]
\draw (0,0.45) -- (0,-0.25)-- (0,0) -- (1.10,0) --(1,0) -- (1,0.45);
\draw (1/6,0.01)--(1/6,-0.01);
\draw (0.25,0.01)--(0.25,-0.01); \draw (0.75,0.01)--(0.75,-0.01);
\draw (0.5,0.01)--(0.5,-0.01);
 \draw (3/8,0.01)--(3/8,-0.01);  \draw (5/8,0.01)--(5/8,-0.01);
\draw (0.99,5/12)--(1.01,5/12);
\draw (1/4,-5/24)--(1,5/12);
\draw (1/6,-2/9)--(1/4,-5/24);
\draw (0,-1/12)--(1/6,-2/9);
\draw (-0.01,-1/6)--(0.01,-1/6);
\draw (-0.01,-0.25)--(0.01,-0.25);\draw (-0.01,-1/12)--(0.01,-1/12);
\draw (0,0) node[left] {$0$} (1/6,0) node[above] {$\frac{3\a+1}{6\a+4}$} (1/4,0) node[below] {$\frac{1}{4}$}  (1/2,0) node[below] {$\frac12$} (3/4,0) node[below] {$\frac{3}{4}$} (1,0) node[below] {$1$} (1.1,0) node[below] {$\frac{1}{p'}$};
\draw (1,5/12) node[right] {$\frac{\alpha+1}{2}$} (0,-0.25) node[left] {\tiny $-1/4$} 
(0,-1/12) node[left] {\tiny $\a/2$}      (0,-1/6) node[left] {\tiny $-1/6$}; 
\end{tikzpicture}
    \caption{$-1/3 < \a <0  $ }
    \label{FigLem2.2}
    \end{minipage}
    \hfill
    \begin{minipage}[b]{0.45\textwidth}
         \begin{tikzpicture}[scale=6]
    \draw (0,0.45) -- (0,-0.25)-- (0,0) -- (1.10,0) --(1,0) -- (1,0.45);

\draw (0.25,0.01)--(0.25,-0.01); 
\draw (0.5,0.01)--(0.5,-0.01);
\draw (1/4,-7/48)--(1,7/24);
\draw (0,-1/8)--(1/4,-7/48);
\draw (-0.01,-1/6)--(0.01,-1/6);
\draw (-0.01,-1/8)--(0.01,-1/8);
\draw (0.99,7/24)--(1.01,7/24);
\draw (1,7/24) node[right] {$\frac{\alpha+1}{2}$};
\draw (0,0) node[left] {$0$} (1/4,0) node[below] {$\frac{1}{4}$}  (1/2,0) node[below] {$\frac12$} (1.10,0) node[below] {$\frac{1}{p'}$};
\draw(0,-1/8) node[left] {\tiny $-(3\alpha+2)/6$}      (0,-1/6) node[right] {\tiny $-1/6$}; 
\end{tikzpicture}
    \caption{$-1/2 \le \alpha \le -1/3 $ }
    \label{FigLem2.2}
    \end{minipage}
\end{figure}

Within the corresponding admissible range of exponents, the radial estimates derived here exhibit improved dependence on the spectral parameter compared to the full Hermite case. This behaviour is consistent with observations made in the work of Imekraz, Robert, and Thomann \cite{IRT16}, where radial Hermite functions are shown to enjoy sharper bounds in the same range of $p$. Below we give the estimates for the $L^p$ norm of the eigenfunctions $\psi^{\a,\b}_k$.

\begin{proposition}\label{PropPsi}
    Let $1\le p\le\infty$. Then for different ranges of $\a,\b\ge-1/2$ we have
     \begin{equation}
    \|\psi^{\a,\b}_k\|_{L^p(d\mu_\b)}\sim k^{\rho(\b,p)},
\end{equation}
where for $\b\ge 0$ and any $\a\ge\b$ 
        \begin{align}\label{Est1a}
     \rho(\b,p):=\begin{cases}
                (\frac{1}{p}-\frac{1}{2})(1+\b),\quad 1\le p\le \frac{4\b+4}{2\b+1},\\
                \frac{\b}{2}-\frac{1}{p}(1+\b),\quad \frac{4\b+4}{2\b+1}\le p\le \infty.
            \end{cases}
        \end{align}
        For $-1/3\le\b\le 0$ and any $\a\ge\b$
        \begin{align}\label{Est1b}
            \rho(\b,p):=\begin{cases}
                (\frac{1}{p}-\frac{1}{2})(1+\b),\quad 1\le p\le 4,\\
                -\frac13+\frac{1}{p}(\frac13+\b)-\frac{\b}{2},\quad 4\le p\le \frac{6\b+4}{3\b+1},\\
                \frac{\b}{2}-\frac{1}{p}(1+\b),\quad \frac{6\b+4}{3\b+1}\le p\le\infty. 
            \end{cases}
        \end{align}
        For $-1/2\le\b\le-1/3$ and any $\a\ge\b$
        \begin{align}\label{Est1c}
            \rho(\b,p):=\begin{cases}
                (\frac{1}{p}-\frac{1}{2})(1+\b),\quad 1\le p\le 4,\\
                -\frac13+\frac{1}{p}(\frac13+\b)-\frac{\b}{2},\quad 4\le p\le \infty. 
            \end{cases}
        \end{align}
        For $\b\ge0$ and $\a<\b$ and all  $p\in [1,p_*)$,
        \begin{align}\label{Est1d}
            \rho(\b,p):=\begin{cases}
                (\frac{1}{p}-\frac{1}{2})(1+\b),\quad 1\le p\le \frac{4\b+4}{2\b+1},\\
                \frac{\b}{2}-\frac{1}{p}(1+\b),\quad \frac{4\b+4}{2\b+1}\le p\le p_*. 
            \end{cases}
        \end{align}
Moreover,
     for $p_*\le p\le\infty$ we have
        \begin{equation*}
            \|\psi^{\a,\b}_k\|_{L^p(d\mu_\b)}=\infty,
        \end{equation*}
        where $p_*=\frac{2\b+2}{\b-\a}.$
\end{proposition}

\begin{proof}
We consider the $p$-norm of $\psi^{\a,\b}_k$ in $L^p(d\mu_\b)$, which is equal to the $p$-norm of $\mathcal{L}^\a_k(x) x^{-\b/2}$ in $L^p(\R^+,x^{\b}dx)$.
Indeed,
 \begin{align}\label{ltol}
     ||\psi_k^{\a,\b}||^p_{L^p(d\mu_\b)}&=\int_0^\infty |\psi_k^{\a,\b}(x)|^p x^{2\b+1}dx=\int_0^\infty |\mathcal{L}^\a_k(x^2) x^{-\b}|^p x^{2\b+1}dx\\
     &=\int_0^\infty |{\mathcal L}_k^{\a}(x) x^{-\b/2}|^p x^\b dx= \int_0^\infty |{\mathcal L}_k^{\a}(x) x^{\b(\frac{1}{p}-\frac12)}|^p dx .\nonumber
 \end{align} 
We want
$$
 \int_0^\infty |{\mathcal L}_k^{\a}(x) x^{\b(\frac{1}{p}-\frac12)}|^p dx=
\int_0^\infty |{\mathcal L}_k^{\alpha'+\beta'} (x) x^{-
\beta'/2}|^p dx
$$
for some $\alpha'$ and $\beta'$ such that $\alpha'+\beta'>-1$.
This gives,
$$
x^{\b(1/p-1/2)}=x^{-\beta'/2}.
$$
Thus, 
$$
\b(1-p/2)=\beta'(-p/2) \quad \mbox{and} \quad \alpha'+\beta'= \a. 
$$
Note that we always have
\begin{equation}\label{Hyp1}
    \alpha'+\beta'= \a >-1,
\end{equation}
and we want $\a'>-2/p$, that is, 
\begin{equation}\label{Hyp2}
 \a'+\frac2p=\a-\b'+\frac2p=\a-\b(1-\frac2p)+\frac2p=(\a-\b)+\frac2p(1+\b)>0.   
\end{equation}
Since $ 1 +\beta \geq 1/2, $  \eqref{Hyp2} holds trivially for all $ \alpha \geq \beta.$ When $ \alpha < \beta$ if we define $ p^\ast =  \frac{2(1+\beta)}{\beta-\alpha} $ then \eqref{Hyp2} holds only for $ 1 \leq p < p^\ast.$ In fact  for $ p \geq p^\ast,$ the functions $ \psi_k^{\alpha,\beta} $ fail to be in $ L^p(\R^+,d\mu_\beta).$ To see this, note that for $0\le x^2 \le 1/\lambda_k$  and any $k$ one has 
$$\psi^{\a,\b}_k(x)\sim  \lambda_k^{\a/2}x^{\a-\b}.$$
Hence, if $p \ge \frac{2(\b+1)}{\b-\a}$, then $ p(\a-\b)+2\b+1 \le -1 $ so that
$$
\|\psi^{\a,\b}_k\|_{L^p(\R^+,d\mu_\b)}^p \ge C\int_0^{1/ \sqrt{\lambda_k}}
x^{p(\alpha-\b)+2\b+1}dx = \infty.
$$
From Lemma \ref{Lag-p-norm} we obtain
\begin{align*}
&(i) \quad \|\psi^{\a,\b}_k \|_{L^p(d\mu_\b)} \sim k^{(1+\b)(\frac1p-\frac12)}
\quad \mbox{if} \quad 
\b(p-2)\le \frac{4-p}{2}
\\
&(ii)\quad  \|\psi^{\a,\b}_k \|_{L^p(d\mu_\b)} \sim k^{\frac\b2-\frac1p(\b+1)}
\quad \mbox{if} \quad 
\b(p-2)> \frac{4-p}{2}
\end{align*}
for p in the interval $1 \le p \le 4$. For $p>4$
\begin{align*}
&(iii) \quad \|\psi^{\a,\b}_k\|_{L^p(d\mu_\b)} \sim k^{\frac{-1}{3}+ \frac{1}{3p}+\b(\frac1p-\frac12)  }
\quad \mbox{if} \quad 
\b\le\frac{4-p}{3p-6}\text{ i.e., } p\le \frac{6\b+4}{3\b+1}
\\
&(iv) \quad \|\psi^{\a,\b}_k\|_{L^p(d\mu_\b)} \sim k^{\frac\b2-\frac1p(\b+1)}
\quad \mbox{if} \quad 
\b>\frac{4-p}{3p-6}\text{ i.e., } p> \frac{6\b+4}{3\b+1}.
\end{align*}
Now we have the $p$-norm of $\psi^{\a,\b}_k$ on $L^p(\R^+,d\mu_\b)$, so we only need to find the ranges of $p$ for which these estimates hold for different ranges of $\b$. 

\begin{itemize}
    \item [(a)]When $\b\ge 0$ we are in cases $(i)$ and $(iv)$, giving exponents \eqref{Est1a}. From $(i)$, we obtain that $1\le p\le\frac{4\b+4}{2\b+1}$ and $(iv)$ gives $p\ge\frac{4\b+4}{2\b+1}$ (since $ \frac{4\b+4}{2\b+1}\ge \frac{6\b+4}{3\b+1}$ whenever $\b\ge0$). 
    \item [(b)]When $-1/3\le\b\le0$, we have the cases $(i),\;(iii)$ and $(iv)$. The ranges of $p$ we get from these cases are $1\le p\le 4$, $4\le p\le \frac{6\b+4}{3\b+1}$ and $p\ge \frac{6\b+4}{3\b+1}$, respectively. 
    \item [(c)]For $-1/2\le\b\le-1/3$, we have cases $(i)$ and $(iii)$ which give the ranges $1\le p \le 4$ and $p\ge 4$, respectively. 
\end{itemize}

\end{proof}

\begin{remark}  In view of the above, when $ \alpha < \beta$ the Bochner-Riesz means of a function $ f \in L^p(\R^+, d\mu_\beta)$ belong to the same space only for  $ 1 \leq p <(2\beta+2)/(\beta-\alpha).$ Hence it is reasonable to assume that $ \alpha \geq \beta $ in our theorems. Otherwise we need to state the case $ -1/2 \leq \alpha < \beta $ each time separately.
\end{remark}

\section{Necessary conditions for Bochner-Riesz summability}
In this section, we verify the necessary conditions for Bochner-Riesz summability for the operator $\Delta_{\a,\b}$. To achieve this goal we define projection operators $P^{\a,\b}_k$ associated to our operator $\Delta_{\a,\b}$. Corresponding to the eigenfunctions $\psi^{\a,\beta}_k$, we define these operators as
\begin{equation}\label{Pdef}
   P^{\a,\beta}_kf:=\langle f,\psi_k^{\a,\b} \rangle_\b \psi_k^{\a,\b},
\end{equation}
acting on spaces $L^p(\R^+,x^{2\b+1}dx)$. Here we set $\langle f,g \rangle_\b=\int_{\R^+}f(x) \bar{g}(x)x^{2\b+1}dx$.

We first describe asymptotics of $L^p-L^2$ norms of the projection operators $P^{\a,\b}_k$ acting on spaces $L^p(\R^+,x^{2\b+1}dx)$ defined by \eqref{Pdef}. 

\begin{proposition}\label{Rmk}
Let $1\le p\le\infty$. Then for different ranges of $\a,\b\ge-1/2$ we have the following $L^p(d\mu_\b)-L^2(d\mu_\b)$ restriction estimate for $P_k^{\a,\b}$. That is, 
\begin{equation*}
  \|P_k^{\alpha,\b}\|_{p \to 2} \sim k^{\rho(\b,p')},
\end{equation*}
where $p'$ is the conjugate exponent to $p$ ($1/p+1/p'=1$).
Note that for $\b\ge0$ and $\a<\b$, $p_*=\frac{2\b+2}{\b-\a}$
so  $p'_*=\frac{2\b+2}{\b+\a+2}$.
    
\end{proposition}

\begin{proof}
    Note that it follows from \eqref{Pdef}
    that $\|P_k^{\a,\b}\|_{p\to 2}= \|\psi_k^{\a,\b}\|_{p'}$. Therefore Proposition \ref{Rmk} is a direct consequence of Proposition \ref{PropPsi}.
 \end{proof}

The following lemma gives the $p$ to $p$ norm of the projection operator $P^{\a,\b}_k$ which is important in obtaining the necessary condition for the Bochner-Riesz means of $\Delta_{\a,\b}$.

\begin{lemma}\label{pn}
Let $1\le p\le\infty$. Then for different ranges of $\b\ge-1/2$ and $\a\ge\b$ we have the following $L^p(d\mu_\b)-L^p(d\mu_\b)$ restriction estimates for $P_k^{\a,\b}$. That is, 
\begin{equation}
    \|P_k^{\alpha,\b}\|_{p \to p} \sim  k^{\tau(\b,p)},
\end{equation}
where
for $\b\ge 0$ 
\begin{equation}\label{pp'-product0}
  \tau(\b,p):=\rho(\b,p)+\rho(\b,p')=    \begin{cases}
    2(\b+1)(\frac1p-\frac12)-\frac{1}{2},  & 
    {{1\le p< \frac{4\b+4}{2\b+3}}},\\
    0, &  {{\frac{4\b+4}{2\b+3}\le p<\frac{4\b+4}{2\b+1}}},\\
   2(\b+1)(\frac{1}{2}-\frac1p)-\frac{1}{2}, & \frac{4\b+4}{2\b+1}\le p\le \infty.
\end{cases}
\end{equation}
For $-1/3<\b<0$
\begin{align}\label{pp'-product1/3}
    \tau(\b,p):=\rho(\b,p)+\rho(\b,p')=
    \begin{cases}
     2(\b+1)(\frac{1}{p}-\frac12)-\frac{1}{2}, &   1\le p<\frac{6\b+4}{3\b+3},\\
    \frac{2}{3}(\frac1p-\frac{1}{2})-\frac16, &  \frac{6\b+4}{3\b+3}\le p< \frac{4}{3},\\
    0, &  \frac{4}{3}<p<4,\\
    \frac{2}{3}(\frac{1}{2}-\frac1p)-\frac16, & 4< p\le \frac{6\b+4}{3\b+1},\\
      2(\b+1)(\frac12-\frac{1}{p})-\frac{1}{2}, & \frac{6\b+4}{3\b+1}<p\le \infty.
\end{cases}
\end{align}
For $-1/2\leq \b\le -1/3$ 
\begin{align}\label{pp'-product1/2}
    \tau(\b,p):=\rho(\b,p)+\rho(\b,p')=
    \begin{cases}
    \frac{2}{3}(\frac1p-\frac{1}{2})-\frac16,
  &  1\leq p\le \frac{4}{3},\\
    0, &  \frac{4}{3}<p<4, \\
     \frac{2}{3}(\frac12 - \frac1p)-\frac16, & 4 \le p \le \infty. 
\end{cases}
\end{align}
Moreover, for $\b\ge0$ and $\a<\b$,
  \begin{equation}\label{Proj3.6}
            \|P_k^{\alpha,\b}\|_{p \to p}=\infty, \text{ for } 1 \le p \le p'_*  \text{ or for } p_* \le p \le \infty.
        \end{equation}
\end{lemma}
\begin{proof}
The proof is a direct consequence of Proposition \ref{PropPsi}
\end{proof}

\begin{theorem}\label{verlp}
 Suppose the  Bochner-Riesz means of order $\delta$ for the operator $\Delta_{\a,\b}$ are uniformly bounded on $L^p
 (\mathbb{R},_+x^{2\b+1}dx)$ for a given $ p, 1\leq p \leq \infty,$ that is
$$
\sup_{R>0}\|\ S_R^{\delta}(\Delta_{\a,\b})\|_{p \to p} \le C < \infty.
$$
Then it necessarily follows that
$$
\delta \ge  \max\left\{0,\quad \frac{2}{3}\left|\frac{1}{2}-\frac{1}{p}\right|-\frac{1}{6}, \quad
2(\b+1)\left|\frac{1}{2}-\frac{1}{p}\right|
-\frac{1}{2}
\right\}.
$$
Moreover, when $\beta\ge 0$ and $\a<\b$, $\sup_{R>0}\|S_R^{\delta}(\Delta_{\a,\b})\|_{p \to p}=\infty$ for $1 \le p \le p'_*  \text{ (or } p_* \le p \le \infty).$
\end{theorem}

\begin{proof}
We start by introducing the distributions $\chi_-^a$, defined by
$$
\chi_-^a=\frac{x_-^a}{\Gamma(a+1)},
$$
where $\Gamma$ is the Gamma function and
$$
x_-^a=|x|^a \quad \mbox{if} \quad x \le 0 \quad \quad \mbox{and}\quad  x_-^a=0 \quad  \mbox{if} \quad x > 0.
$$
The $x_-^a$ are clearly distributions for $\Re a > -1$, and we have
for $\Re a > 0$,
\begin{equation*}
	\frac{d}{dx} x_-^a = -a x_-^{a-1} \implies \frac{d}{dx} \chi_-^a=-\chi_-^{a-1}
\end{equation*}
which we use to extend the family of functions $ \chi_-^a$ to a family of distributions on $\RR$ defined for all $a\in \CC$; see \cite{Ho1} for details.  Since $1-\chi_-^0(x) $ is the Heaviside function, it follows that
\begin{equation}
	\chi_-^{-m} =(-1)^m\delta_0^{(m-1)}, \quad  m=1,2,\ldots ,
	\label{Heaviside}
\end{equation}
where $\delta_0$ is  the $\delta$-Dirac measure. Note that it should not be confused with the order $\delta$ in the statement of the theorem.

Motivated by the above equality we define the distribution $\delta_-^{\nu}$ for  all real exponents $\nu \in \R$ by
the formula
$$
\delta_-^{\nu}=\chi_-^{-\nu-1}.
$$
A straightforward computation shows that for all $w,z\in \CC$
\begin{equation*}\label{aa2}
	\chi_-^w * \chi_-^z=\chi_-^{w+z+1},
\end{equation*}
where $\chi_-^w* \chi_-^z$ is the convolution of the distributions $\chi_-^w$ and $\chi_-^z$
see \cite[(3.4.10)]{Ho1}. It follows from the above relation that
$$
\delta_-^{\nu}*\delta_-^{\mu}=\delta_-^{\nu+\mu}.
$$
Now if $\supp F \subset [0,\infty)$ we define the Weyl fractional
derivative of $F$ of order $\nu$ by the formula
$$
F^{(\nu)}=F*\delta_-^{\nu}
$$
and we note that
$$
F^{(\nu)}*\delta_-^{-\nu}=F*\delta_-^{\nu}*\delta_-^{-\nu}=F,
$$
see \cite[Page 308]{GP}. 
Thus,
$$
F(\Delta_{\a,\b})=\frac{1}{\Gamma(\nu)}\int_0^\infty F^{(\nu)}(t)
(t-{\Delta_{\a,\b}})_+^{\nu-1}dt= \frac{1}{\Gamma(\nu)}\int_0^\infty
F^{(\nu)}(t)t^{\nu-1}S_t^{\nu-1}(\Delta_{\a,\b})dt,\ \ \forall{\nu >0}
$$
for all $F$ supported in the positive half-line. 
Hence  if $\supp F \subset [0,\infty)$ then
\begin{eqnarray}
	\|F(\Delta_{\a,\b})\|_{p\to p}\leq  C\sup_{N}\| S_N^{\nu}(\Delta_{\a,\b})\|_{p\to
		p} \int_0^\infty \left|F^{(\nu+1)}(t)\right|t^{\nu}dt.
	\label{e4.1}
\end{eqnarray}

Consider now a function $\eta \in C_c^\infty(\R)$  such that $\eta(0)=1$ and  $\supp\,
\eta \subset [-{1},{1}]$. We define a sequence of functions
$F_k$ by the formula
\begin{equation*}
    F_k(\lambda)=\eta((\lambda-\lambda_k)).
\end{equation*}
Observe that
\begin{equation}\label{F_k}
\int_0^\infty \left|F_k^{(\nu+1)}(t)\right|t^{\nu}dt\sim\lambda_k^\nu
\sim k^{\nu},  
\end{equation}
and  $F_k(\lambda_l) =1$, if $k=l$ and
$F_k(\lambda_l) =0$ otherwise.
Hence
$$
F_k(\Delta_{\a,\b})f=\sum_{l=1}^{\infty}F_k(\lambda_l)\langle{\psi}_l^{\alpha,\b},f\rangle_\b{\psi}_l^{\alpha,\b}=\langle{\psi}_k^{\alpha,\b},f\rangle_\b{\psi}_k^{\alpha,\b}=P^{\a,\b}_kf.
$$
Thus, by Lemma~\ref{pn} and estimates \eqref{e4.1} and \eqref{F_k}  
we have
\begin{equation*}
    \|P^{\a,\b}_k\|_{p\to p}\sim k^{\tau(\beta,p)}\lesssim k^\nu.
\end{equation*}
In other words, if the Bochner-Riesz means of order $\nu$ are uniformly bounded in \eqref{e4.1} then $\nu\ge\tau(\b,p)$, which is equivalent to saying
\begin{equation*}
    \nu\ge\max\left\{0,\quad \frac{2}{3}\left|\frac{1}{2}-\frac{1}{p}\right|-\frac{1}{6}, \quad
2(\b+1)\left|\frac{1}{2}-\frac{1}{p}\right|
-\frac{1}{2}
\right\}.
\end{equation*}
Note that when $\b\ge 0$ and $\a<\b$, we have $\|P^{\a,\b}_k\|_{p\to p}=\infty$ for $1\le p\le p_*'$ . So, by \eqref{e4.1} $\sup_{R>0}\|S^\delta_R(\Delta_{\a,\b})\|_{p\to p}=\infty$ for $1\le p\le p_*'$.

\end{proof}

Below we give systematic diagrams to describe the necessary condition (Theorem~\ref{verlp}) for the boundedness of the Bochner-Riesz means of order $\delta$.

For $\b \ge 0 $, the region $A$ in Figure~\ref{fig4} is the region of uniform boundedness for the Bochner-Riesz means, and we have unboundedness in regions 
\begin{equation}\label{b}
    B=\left\{(p,\delta)\colon \, \delta \le \max \Big\{0,
2(\b+1)\big|\frac{1}{2}-\frac{1}{p}\big|-\frac{1}{2} \Big\}\right\}
\end{equation}
and 
\begin{equation}\label{c}
    C=\left\{(p,\delta)\colon \, \delta \le \max \Big\{0,
\frac{2}{3}\big|\frac{1}{2}-\frac{1}{p}\big|-\frac{1}{6} \Big\}\right\}.
\end{equation}

\begin{figure}[h!]
    \centering
    \begin{minipage}[b]{0.45\textwidth}
\begin{tikzpicture}[scale=6]
\draw (0,0) -- (1,0) -- (1,0.6) -- (0,0.6) -- cycle;
\filldraw[fill=black!40] (0,0) -- (0,1/6) -- (1/4,0) -- cycle;
\filldraw[fill=black!40] (1,0) -- (1,1/6) -- (3/4,0) -- cycle;
\draw[pattern=custom north west lines]  (0,0) -- (0,1/2) -- (2/5,0) -- cycle;
\draw[pattern=custom north west lines] (1,0)--(1,1/2)--(3/5,0)--cycle;
\draw (0,0) node[below] {$0$} (1/4,0) node[below] {$\frac{1}{4}$} (2/5,0) node[below] {$\frac{2\b+1}{4\b+4}$}
(3/5,0) node[below] {$\frac{2\b+3}{4\b+4}$} (3/4,0) node[below] {$\frac{3}{4}$} (1,0) node[below] {$1$} (1.1,0) node[below] {$\frac1p$};
\draw (0,1/6) node[left] {$\frac{1}{6}$} (0,1/2) node[left] {\small $\b+\frac{1}{2}$} (0,0.6) node[left] {$\delta$}; 
\draw (0.5,0.45) node {$A$} (0.2,0.2) node {$B$} (0.05,0.05) node {$C$} (0.8,0.2) node {$B$} (0.95,0.05) node {$C$};
\end{tikzpicture}
    \caption{For $\b\ge0$}
    \label{fig4}           
    \end{minipage}
    \hfill
    \begin{minipage}[b]{0.45\textwidth}
        \begin{tikzpicture}[scale=6]
\draw (0,0) -- (1,0) -- (1,0.6) -- (0,0.6) -- cycle;
\filldraw[fill=black!40] (0,0) -- (0,1/3) -- (1/2,0) -- cycle;
\draw[pattern=custom north west lines] (0,0) -- (0,1/2) -- (1/3,0);
\draw (0,0) node[below] {$0$} (1/5,0) node[below] {$\frac{3\b+1}{6\b+4}$} (1/2,0) node[below] {$\frac{1}{4}$}  (1,0) node[below] {$\frac12$} (1.1,0) node[below] {$\frac1p$};
\draw (0,1/3) node[left] {$\frac{1}{6}$} (0,1/2) node[left] {\small $\b+\frac{1}{2}$} (0,0.6) node[left] {$\delta$}; 
\draw (0.5,0.45) node {$A$} (0.05,0.35) node {$B$} (0.05,0.05) node {$C$} ;
\draw[dashed] (1/5,1/5)--(1/5,0);
\end{tikzpicture}
  \caption{For $-1/3<\b<0$}
    \label{fig5}
    \end{minipage}
\end{figure}

In Figure~\ref{fig5} we focus on half of the region because of symmetry around $p=2$. Therefore, for $-1/3 < \b< 0 $, we have the dotted region below that goes down to $\delta=\frac{1}{6}$. We show that in regions $B$ and $C$ we again have unboundedness of the Bochner-Riesz means and boundedness in the new region $A$.

For $-1/2\le \b \le -1/3$ one has $B \subset C$.

\bigskip
\subsection{Laguerre expansions of Hermite type}
Recall that $ \tilde{\Delta}_{\alpha}=\Delta_{\a,-1/2}$ is the Laguerre operator of Hermite type whose eigenfunctions are provided by $ h_k^\alpha$ known as  Laguerre functions of Hermite type.
\begin{proposition}\label{pr2.2}
Let $1\le p\le\infty$ and $\a\ge-1/2$. Then the system determined by $h_k^\a$ satisfies the following asymptotics for $L^p(dx)$ norms 
\begin{align*}
            \|h_k^\a\|_{L^p(dx)}\sim k^{\rho(-1/2,p)}=\begin{cases}
            k^{\frac{1}{2p}-\frac14},
            \quad 1\le p\le 4\\
              k^{-\frac{1}{6p}-\frac{1}{12}}  ,\quad 4\le p\le \infty.
            \end{cases}
        \end{align*}
\end{proposition}
\begin{proof}

Since $ h_k^\alpha(x)=\mathcal{L}_k^\alpha(x^2) x^{1/2} = \psi^{\a,-1/2}_k(x) $, we have
\begin{equation*}
    \|h^\a_k\|_p=\|\psi^{\a,-1/2}_k\|_p\sim k^{\rho(-1/2,p)}.
\end{equation*}
\end{proof}

Now corresponding to eigenfunctions $h^\a_k$, we define the projection operators $\widetilde{P}^\alpha_k$ in the following way
\begin{equation*}
    \widetilde{P}^\a_kf=P^{\a,-\frac12}_kf=\langle f,\psi_k^{\a,-\frac12}\rangle \psi_k^{\a,-\frac12}=\langle f,h^\a_k\rangle h^\a_k,
\end{equation*}
where $\langle f,g \rangle=\int_{\R^+}f(x) \bar{g}(x)dx$. Then in the following lemma we give the estimates on the $p$ to $p$ norm of the operator $\widetilde{P}^\a_k$.

\begin{lemma}
[$p$ to $p$ norm of the Projection operator]
Let $1\le p\le \infty$. For $\a\ge-1/2$
\begin{equation}\label{pp'h}
  \|\widetilde{P}^\a_k\|_{p \to p}=  \| h_k^\alpha\|_p \|h_k^\alpha\|_{p'}\sim
    \begin{cases}
      k^{-1/2+2/(3p)}, & 
    {{1\le p< \frac{4}{3}}},\\
    1, &  {{\frac{4}{3}< p<4}},\\
     k^{1/6-2/(3p)}, &  4< p \le \infty .
\end{cases}
\end{equation}

\end{lemma}
\begin{proof}
The proof is a direct consequence of Proposition \ref{pr2.2}
\end{proof}

\begin{proposition}\label{verlp1}
 Let $1 \le p \le \infty$. If the  Bochner-Riesz means of order $\delta$ for the operator  $\tilde{\Delta}_{\alpha}$ are uniformly bounded on $L^p$, that is
$$
\sup_{R}\| S_R^{\delta}(\tilde{\Delta}_{\alpha})\|_{p \to p} \le C < \infty,
$$
then it necessarily follows that
$$
\delta \ge  \max\left\{0,\quad \frac{2}{3}\left|\frac{1}{2}-\frac{1}{p}\right|-\frac{1}{6} 
\right\}
.$$

\end{proposition}

\begin{proof}
The proof is a straightforward consequence of Theorem~\ref{verlp} by taking $\b=-1/2$ there and using the estimates from Lemma 3.2 for the corresponding projection operator.
\end{proof}

\subsection{Laguerre expansions of convolution type} As we have already mentioned expansions in terms of $ \varphi_k^\alpha(x) = \psi_k^{\a,\a}(x) $ are known as Laguerre expansions of convolution type for the following reason. For functions $ f \in L^p(\R^+, d\mu_\a)$ and $ g \in L^1(\R^+,d\mu_\alpha),$ Markett \cite{mark2}  defined a convolution structure $ f \ast_\a g$ for any $ \a \geq -1/2.$ The convolution was defined in terms of a Laguerre translation operator $ T_x^\a$ which was motivated by the product formula for Laguerre polynomials. Thus we have
$$ f \ast_\a g(x) = \int_0^\infty f(y)\, T_x^\a g(y)\, d\mu_\a .$$
The Laguerre translation operator has the property that 
$$T_x^\a \psi_k^\a(y)  =\frac{\Gamma(k+1)\Gamma(\a+1)}{\Gamma(k+\a+1)} \psi_k^\a(x)\, \psi_k^\a(y) ,\,\,\,\, \psi_k^\a(x) =  L_k^\a(x^2)\, e^{-\frac{1}{2}x^2}.$$
In view of this the expansion of a function $ f \in L^p(\R^+, d\mu_\a)$ in terms of $ \varphi_k^\a(x) $ can be written in the form
$$ f(x)= \sum_{k=0}^\infty \langle f, \varphi_k^\a \rangle \, \varphi_k^\a(x) = \sum_{k=0}^\infty  \, f \ast_\a \psi_k^\a(x).$$ 
Therefore, the Ces{\`a}ro means of order $ \delta $ associated to the above expansion takes the form
$$ \sigma_N^\delta f(x) = \frac{1}{A_N^\delta} \sum_{k=0}^N A_{N-k}^\delta f\ast_\a \psi_k^\a(x) = f \ast_\a K_N^\delta $$
where the kernel $ K_N^\delta $ can be calculated explicitly. Indeed,
$$ K_N^\delta(x) = \frac{1}{A_N^\delta} \sum_{k=0}^N A_{N-k}^\delta \, \psi_k^\a(x) =  \frac{1}{A_N^\delta} L_N^{\a+\delta+1}(x^2)\, e^{-\frac{1}{2}x^2}.$$ Under the assumption that $ \alpha \geq 0,$ Markett \cite{mark2} proved that the Laguerre translation operator $ T_x^\a $ is bounded on $ L^1(\R^+,d\mu_\a).$ Consequently, it follows that
$$ \| \sigma_N^\delta f\|_{L^1(d\mu_\a)} \leq \|K_N^\delta\|_{L^1(d\mu_\a)}\, \|f\|_{L^1(d\mu_\a)}.$$ 
He then used the fact that  $\|K_N^\delta\|_{L^1(d\mu_\a)}$ is uniformly bounded for $ \delta > \a+1/2 $ under the assumption that $ \a \geq 0$ to complete his proof.\\

However, as can be easily read from Proposition \ref{PropPsi} the uniform bound for $\|K_N^\delta\|_{L^1(d\mu_\a)}$ is valid for any $ \a \geq -1/2 $ and $ \delta> \a+1/2.$  In view of the necessary conditions proved in Theorem \ref{verlp} this simply means that the Laguerre translation operator $ T_x^\a $ is not bounded on $ L^1(\R^+, d\mu_\a) $ for $ -1/2 \leq \a \leq -1/3.$  Although it is reasonable to expect boundedness of $ T_x^\a $ for $ -1/3 < \a <0 $, this remains an open problem. A positive answer to this problem will give a simple proof of Theorem \ref{Th:Sufficient}   for the case $ \a=\beta > -1/3.$

\section{Convergence of the partial sums}

In this section we study the norm convergence of the  partial sums defined by
$$ \sigma_N(\Delta_{\alpha,\beta})f(x) = \sum_{k=0}^N \langle f, \psi_k^{\alpha,\beta} \rangle_\beta\, \psi_k^{\alpha,\beta}(x),\,\,\, f \in  L^p(\R^+, d\mu_\beta).$$
This definition coincides with the definition of Bochner-Riesz means of order 0 for the operator $\Delta_{\a,\b}$.
For $ 1 \leq p <\infty,$ the norm convergence follows  from the uniform boundedness of the partial sum operators. The range of  $ p $ for which  $  \sigma_N(\Delta_{\alpha,\beta})$ are uniformly bounded on 
$ L^p(\R^+, d\mu_\beta)$ depends on the pair $(\alpha,\beta).$  In this section, we give two proofs for the main theorem. The first is a direct consequence of Theorem~\ref{Mucken} by Muckenhoupt and the second proof gives a more refined version of the same, which involves techniques from Herz's work \cite{Herz} where he treated the  problem of inverting the Hankel transforms.\\

As we plan to use the following theorem of Muckenhoupt, let us fix some notation. Just for this section let us use the notation
$$ S_N g(x) = \sum_{k=0}^N (g, \mathcal{L}_k^\alpha)\, \mathcal{L}_k^\alpha(x),\,\,  g \in L^p(\R^+, dx)$$
for the partial sums associated to standard Laguerre expansions. Note that  $ S_N g$ is the same as $ \sigma_N(\Delta_{\a,0})f $ if we define $ f(x) = g(x^2).$ We will use $ \sigma_N(\Delta_{\a,\b})f$ for the partial sums associated to the expansions in terms of $ \psi_k^{\a,\b}.$
We now recall the following theorem of Muckenhoupt  which is a simplified version of \cite[Theorem 7]{Muck2}.

\begin{theorem}[Muckenhoupt'70]\label{Mucken}
 Let $1<p<\infty$ ($p\ne\frac43$ and $p\ne 4$), $\eta>-1$ and $u(x)=x^a$.
 Assume that
 \begin{equation*}
   -\frac{1}{p}+\max\Big(-\frac{\eta}{2},\frac14\Big)<a<1-\frac{1}{p}+\max\Big(-\frac{\eta}{2},\frac14\Big).
 \end{equation*}
Then there exists a constant $C$, independent of $f$ and $N$ such that 
 \begin{equation*}
     \int_0^\infty |S_N f(x)u(x)|^pdx\le C\int_0^\infty|f(x)u(x)|^p dx.
 \end{equation*}
\end{theorem}

Following is one of our main theorems which gives the Bochner-Riesz means of order 0 for $\Delta_{\a,\b}$ and the proof follows from Theorem~\ref{Mucken}.
\begin{theorem}\label{MainThm}
	Suppose that $-1/2 \le \b< 0$.
	Then $$ \sup_{N} \|\sigma_N(\Delta_{\a,\b})f\|_p\le C\|f\|_p
	$$ in $L^p(\R^+,d\mu_\b)$ for all $4/3<p<4$. 
	Next, suppose that $\b\ge 0 $.
	Then 
	$$ \sup_{N} \|\sigma_N(\Delta_{\a,\b})f\|_p\le C\|f\|_p
	$$  in $L^p(\R^+,d\mu_\b)$
	for all $p$ such that $|\frac{1}{p}-\frac12| < \frac{1}{4(\b +1)}$.
\end{theorem}
\begin{proof}

Note that by defining $\tilde{\psi}^{\a,\b}_k(x)=\psi^{\a,\b}_k(\sqrt{x})$ we have
\begin{equation*}
  \tilde{\psi}^{\a,\b}_k(x)=\psi^{\a,\b}_k(\sqrt{x})  =\LL^\a_k(x)x^{-\b/2}.
\end{equation*}
It is clear that the expansions in terms of $\psi^{\a,\b}_k$ in $L^p(\R^+,x^{2\b+1}dx)$ are equivalent to expansions in terms of $\tilde{\psi}^{\a,\b}_k$ in $L^p(\R^+,x^{\b}dx)$. Thus, the uniform boundedness of $\sigma_N(\Delta_{\a,\b})$ on $L^p(\R^+,x^{2\b+1}dx)$ is equivalent to the uniform boundedness of $\sigma_N$ on $L^p(\R^+,x^\b dx)$ where
	\begin{equation*}
		\sigma_N f(x)=\sum_{k=0}^N( f,\tilde{\psi}^{\a,\b}_k)\, \tilde{\psi}^{\a,\b}_k(x),\quad ( f,\tilde{\psi}^{\a,\b}_k)=\int_0^\infty f(x)\,\tilde{\psi}^{\a,\b}_k x^{\b}dx.
	\end{equation*}
Since \begin{equation*}
    \int_0^\infty f(x)\,\tilde{\psi}^{\a,\b}_k(x)x^{\b}dx= \int_0^\infty f(x)\LL^\a_k(x)x^{\b/2}dx
\end{equation*}
we have the relation,
\begin{equation*}
   \sigma_N f(x)=x^{-\b/2}S_{N} g(x), \text{ where } g(x)=x^{\b/2}f(x).
\end{equation*}
Equivalently,
\begin{equation*}
   \sigma_Nf(x)=x^{-\b/2}S_Nx^{\b/2}f(x).
\end{equation*}
For $\b\ge{-1}/2$, we want to find the range of $p$ for which 
\begin{equation}
    \int_0^\infty|\sigma_Nf(x)|^p x^{\b}dx\le C\int_0^\infty|f(x)|^p x^\b dx
\end{equation}
which is equivalent to the weighted norm inequality
\begin{equation}\label{WLE}
    \int_0^\infty|S_Ng(x)|^p x^{\b(1-p/2)}dx\le C\int_0^\infty|g(x)|^p x^{\b(1-p/2)} dx.
\end{equation}
The case $ \b = 0 $  follows from the result of Askey-Wainger. For  $ \b \neq 0,$ we can deduce the above  inequality from Muckenhoupt's theorem for partial sums associated to standard Laguerre expansions stated above.
By taking $\eta=\beta\ge-1/2$ in Theorem~ \ref{Mucken} , the choice
	\begin{equation*}
		-\frac{1}{p}+\frac14 < a <-\frac{1}{p}+\frac34
	\end{equation*}
	is allowed. 
    Finally, we want to choose $a=\b(1/p-1/2)$, since $u(x)=x^a$ in Theorem~ \ref{Mucken}.
    Now, this choice is possible whenever 
    \begin{equation*}
        -\frac{1}{p}+\frac14<\b(\frac{1}{p}-\frac12)<-\frac{1}{p}+\frac34, 
    \end{equation*} that is, whenever
    \begin{equation*}
        \frac{4(\b+1)}{2\b+3}<p<\frac{4(\b+1)}{2\b+1}.
    \end{equation*}
    This holds true whenever $\b\ge-1/2$, which is our assumption. \\Since we have 
    \begin{equation*}
        \frac{4(\b+1)}{2\b+3}<\frac{4}{3}<p<4<\frac{4(\b+1)}{2\b+1}, \text{ for } -\frac{1}{2}\le\b<0,
    \end{equation*}
    and 
    \begin{equation*}
        \frac{4}{3}<\frac{4(\b+1)}{2\b+3}<p<\frac{4(\b+1)}{2\b+1}<4, \text{ for } \b>0
    \end{equation*}
    with our choice of $a$ we get the weighted norm inequality for Laguerre polynomial expansion. That is the same as the weighted norm inequality in \eqref{WLE}.
\end{proof}

When $\a=\b$, we get the following result as a consequence of  Theorem~\ref{MainThm}  for the convergence of partial sums of $\sigma_N(\Delta_{\alpha,\alpha})$.
 
 \begin{coro}\label{partial} Assume that $ \alpha \geq -1/2.$  When $ \alpha\geq 0,$  the inequality \eqref{WLE} (with $\b=\a$) is valid and  hence the partial sum operators $ \sigma_N(\Delta_{\alpha,\alpha})$ associated to expansions  in terms of $ \varphi_k^\alpha $ are uniformly bounded on $ L^p(\R^+,d\mu_\alpha)$ for all $ 4(\alpha+1)/(2\alpha+3) < p < 4(\alpha+1)/(2\alpha+1).$ When $ -1/2 \leq \alpha \leq 0,$ the same conclusion is valid for all $ 4/3 < p <4.$
 \end{coro}

 \begin{remark} We can also show that the conditions in the above corollary are necessary. The uniform boundedness of  $ \sigma_N(\Delta_{\alpha,\alpha})$ implies that
 $$  | \langle f, \varphi_k^\alpha \rangle_\alpha| \, \| \varphi_k^\alpha \|_{L^p(\R^+,d\mu_\alpha)} \leq  C \|f\|_{L^p(\R^+,\,d\mu_\alpha)} ,$$
and hence, the boundedness of the following is necessary:
 $$\| \varphi_k^\alpha \|_{L^p(\R^+,\,d\mu_\alpha)}\| \varphi_k^\alpha \|_{L^{p^\prime}(\R^+,\,d\mu_\alpha)}  \leq  C.$$
 The conditions on $ p $ can be read from the upper and lower bounds for various  norms of $ \varphi_k^\alpha $ stated in  Proposition \eqref{basiclemma}.

 The above calculation provides the essential intuition to establish the necessary condition stated in Theorem~\ref{verlp}.
 \end{remark}

 \subsection{A splitting lemma and some preparatory results} \label{Subsec:SplitLem} The kernel of the partial sum operator $ \sigma_N(\Delta_{\alpha,\beta}) $ involves  functions of the form $ \psi_N^{\alpha,\beta}(x)\, \psi_{N+1}^{\alpha,\beta}(y) $
, and hence their behaviour in various regions depends on the asymptotic properties of the Laguerre functions.  From the estimates in \eqref{psi} it is clear that  the Laguerre function $ \psi_N^{\alpha,\beta}(x) $ has four distinct behaviour depending on whether $ x^2 \in I_j, j=1,2,3,4 $ where $ I_j $ are the intervals  $ [0,N/4], [N/4, N/2], [N/2, 2N] $ and $ [2N, \infty).$  We let $ \chi_N^j(x)  = \chi_{I_j}(x^2) $ where $ \chi_{I_j} $ stand for the characteristic function of the interval $ I_j.$  The following lemmas are rather elementary but play a significant role in our study.

\begin{lemma}\label{elem1}
   Suppose that $\mu$ is a positive measure on $\R^+$ and $N>0$.   For  a linear operator $ T $ acting on $L^p(\R^+, d\mu)$  let $ \|T\|_{p\to p} $ be its operator norm. Then the following inequality holds:
\begin{align*}
    \big\|\chi_N^i T\chi_N^j\big\|_{p\to p} \le \big\|T\big\|_{p\to p}
    \le & 4\big\|(\chi_N^{1}+\chi_N^2)T(\chi_N^{1}+\chi^2_{N}) \big\|_{p\to p}  + 4 \,\big\|(\chi_N^{2}+\chi_N^3)T(\chi_N^2+\chi_N^3) \big\|_{p\to p}  \\
    &+\big\|\chi_N^1T\chi_N^3\big\|_{p\to p} + \big\|\chi_N^3T\chi_N^1\big\|_{p\to p}+
    4\big\|T\chi_N^4\big\|_{p\to p} + 4\big\|\chi_N^4T\big\|_{p\to p}
\end{align*}  
for all  $1 \le p \le \infty$ and $1 \le i,j \le 4$.
\end{lemma}

For any bounded function $ F $ defined on the spectrum of $ \Delta_{\alpha,\beta} $ we define $ F(\Delta_{\alpha,\beta}) $ by the spectral theorem, namely
$$ F(\Delta_{\alpha,\beta})f(x)= \sum_{k=0}^\infty  F(\lambda_k)\, \langle f, \psi_k^{\alpha,\beta}\rangle_\beta\, \psi_k^{\alpha,\beta}(x),\,\,\, \lambda_k=4k+2\a+2.$$

\begin{lemma}\label{2N} Let $ F $ be a bounded function such that $\supp F \subset [0,4N/3]$. Then for any  $1 \le p \le \infty$ there exists a constant $ C $ independent of $ N $ such that
$$
\|F(\Delta_{\a,\b})\chi_N^4\|_{{L^p(d\mu_\b)}\to {L^p(d\mu_\b)} } \leq
C\|F\|_{\infty} \quad \mbox{and} \quad \|\chi_N^4F(\Delta_{\a,\b})\|_{{L^p(d\mu_\b)}\to {L^p(d\mu_\b)} }
\leq
C\|F\|_{\infty}$$

\end{lemma}

\begin{proof}
We only verify the first inequality. The second is a straightforward modification of the first.
By the definition of spectral multipliers
$$
F({\Delta_{\a,\b}})\chi_N^4f = \sum_{1\le \lambda_k \le 4N/3}
F(\lambda_k)\, \langle \chi_N^4f,
\psi^{\a,\b}_k \rangle_\b\, {\psi^{\a,\b}_k},
$$
and hence it follows that
\begin{eqnarray*}
\|F(\Delta_{\a,\b})\chi_N^4\|_{p \to p} \leq\sum_{1\le \lambda_k \le 4N/3}
|F(\lambda_k)|  \,\|{\psi^{\a,\b}_k}\|_{p}\, 
\|\chi_N^4\psi^{\a,\b}_k\|_{p'}.
\end{eqnarray*}
By Proposition~\ref{PropPsi}, we know that $
\|\psi^{\a,\b}_k\|_{p} \leq Ck^{\rho(\b,p)}
$
and since $\supp F \subset [0, 4N/3],$  in the sum above it is enough to  consider only  such $k$ so that
$3\lambda_k/2 \leq 2N $. Hence by \eqref{psi}
$$
\|\chi_N^4\psi^{\a,\b}_k\|_{p'} ^{p'} \le C \int_{2N}^\infty e^{-p'\gamma x} dx \leq C'\exp(-p'\gamma N)
$$
and consequently, we have the estimate
$$
\|F(\Delta_{\a,\b})\chi_N^4\|_{p \to p}
\leq CN^{\rho(\b,p)}\exp(-\gamma N)\Bigg(\sum_{\lambda_k \le 4N/3} |F(\lambda_k)|\Bigg).
$$
Since the number of eigenvalues below $N$ is of order $N,$  we have
$$
\sum_{\lambda_k \le 4N/3} |F(\lambda_k)| \leq CN\|F\|_{\infty}.
$$
This proves Lemma~\ref{2N}.
\end{proof}

\subsection{An alternate proof for the convergence of partial sums} In this subsection we present an alternate proof of Theorem \ref{MainThm} following the ideas of Herz.  Let us fix $ \alpha $ and $ \beta $ and write $ \sigma_N $ instead of  $ \sigma_N(\Delta_{\alpha,\beta}) $ throughout this subsection. The partial sum operator $ \sigma_N $ is an integral operator:
$$ \sigma_Nf(x) = \int_0^\infty \sigma_N(x,y)\, f(y)\, y^{2\beta+1}\,dy,\,\,\,\,\,   \sigma_N(x,y) = \sum_{k=0}^N \psi_k^{\alpha,\beta}(x)\, \psi_k^{\alpha,\beta}(y).$$
As done by the earlier authors (Askey-Wainger and Muckenhoupt) we also use Pollard's idea of using the expression 
$$ \frac{1}{2}(\sigma_N(x,y) + \sigma_{N+1}(x,y)) = \sigma_N(x,y)+ e_N(x,y), \,\,\, e_N(x,y) =\frac{1}{2}  \psi_{N+1}^{\alpha,\beta}(x)\, \psi_{N+1}^{\alpha,\beta}(y). $$  
The uniform boundedness of the operators  with kernel $ e_N(x,y) $ holds if and only if 
$$\| \psi_k^{\alpha,\beta} \|_{L^p(\R^+,\,d\mu_\beta)}\| \psi_k^{\alpha,\beta} \|_{L^{p^\prime}(\R^+,\,d\mu_\beta)}  \leq  C,$$ 
which is true for $\beta\ge 0 $ whenever $\left|\frac{1}{p}-\frac12\right|< \frac{1}{4\b+4}$. For $-\frac{1}{2}\le\b<0$, it holds whenever $\frac43<p<4$.
Therefore, in view of the estimates stated in Proposition \ref{PropPsi} and Lemma~\ref{pn}, we only need to  consider the operator 
$$ D_Nf(x) = \int_0^\infty D_N(x,y) \, f(y)\, y^{2\beta+1}\, dy,\,\,\,   D_N(x,y) = \frac{1}{2} ( \sigma_N(x,y) + \sigma_{N+1}(x,y)).$$
Recall that the Laguerre polynomials $ L_k^\alpha $ satisfy the identity 
$$ \sum_{k=0}^N \frac{\Gamma(k+1)}{\Gamma(k+\alpha+1)}\, L_k^\alpha(x) \, L_k^\alpha(y)  =\frac{\Gamma(N+2)}{\Gamma(N+\a+1)}\frac{1}{x-y}
\big( {L}^\a_N(x){L}^\a_{N+1}(y) - {L}^\a_{N+1}(x){L}^\a_N(y) \big) .$$
In view of this, we have
$$\sigma_N(x,y)  =\big( (N+1)(N+\a+1)  \big)^{1/2}\frac{1}{x^2-y^2}\big( \psi^{\a,\b}_N(x)\psi^{\a,\b}_{N+1}(y) - \psi^{\a,\b}_{N+1}(x)\psi^{\a,\b}_N(y) \big).$$
It follows that we can write the kernel as a sum, $ D_N(x,y) = \sum_{j=1}^3 D_{N,j}(x,y) $, where
$$ D_{N,1}(x,y) =\frac{c_{N+1}}{2(x^2-y^2)}\psi^{\a,\b}_{N+1}(x)[\psi^{\a,\b}_{N+2}(y)-\psi^{\a,\b}_N(y)],$$
$$ D_{N,2}(x,y) = \frac{c_{N+1}}{2(x^2-y^2)}\psi^{\a,\b}_{N+1}(y)[\psi^{\a,\b}_N(x)-\psi^{\a,\b}_{N+2}(x)], $$
$$ D_{N,3}(x,y) = \frac{(c_{N+1}-c_N)}{2(x^2-y^2)}\big[\psi^{\a,\b}_{N+1}(x)\psi^{\a,\b}_N(y)-\psi^{\a,\b}_{N+1}(y)\psi^{\a,\b}_N(x)\big].$$
In the above  $c_N=((N+1)(N+\a+1))^{1/2}$ so that $ c_N+c_{N+1}= O(N)$ and $ ( c_{N+1}-c_N)=O(1)$. \\

Let $ D_{N,j} $ stand for the operator with kernel $ D_{N,j}(x,y).$ To study the uniform boundedness of these operators, we use the splitting of Lemma~\ref{elem1}. We divide each of our operators into four different pieces and study the boundedness of each  piece using the right hand side of the splitting  Lemma~\ref{elem1}. \\

Note that the $L^p$ boundedness of operators $\chi^4_ND_N$ and $D_N\chi^4_N$ comes from Lemma~\ref{2N} and the proof for the remaining part of the splitting lemma for the operator $D_N$ is given below in Lemma~\ref{1to3} and Theorem~\ref{th10.a}.

\begin{lemma}\label{1to3} 
Suppose $\a,\b \ge -1/2$.  
Then for any  $1 \le p \le \infty$
$$
\|\chi_N^1   D_N\chi^3_N\|_{{L^p(d\mu_\b)}\to {L^p(d\mu_\b)} } +\|\chi^3_N D_N\chi^1_N\|_{{L^p(d\mu_\b)}\to {L^p(d\mu_\b)} } 
\leq C N^{\tau(\b,p)}
$$
for all $N >0$, where $\tau(\b,p)$ is defined in Lemma~\ref{pn}. 
\end{lemma}

\begin{proof}
The argument to obtain the estimates for $\|\chi^1_N   D_N\chi^3_N\|_{p\to p}$ and $\|\chi^3_N   D_N\chi^1_N\|_{p\to p}$ follows in the same way, therefore we consider only the former operator.
We know that $D_N= \sum_{j=1}^3D_{N,j}$. Therefore,
for $\|\chi^1_N D_{N,1}\chi^3_N\|_{p\to p}$ note that if  $ y\in \supp  \,\chi^3_N $ and
 $x \in \supp  \,\chi^1_N $ then $\frac{1}{|x^2-y^2|}\sim\frac{1}{N}$, which cancels the growth in $N$ coming from $c_{N+1}$. This implies,
 \begin{equation*}
    \|\chi^1_N D_{N,1}\chi^3_N\|_{p\to p}\lesssim \|\chi_N^1\psi^{\a,\b}_{N+1}\|_p \;\|(\psi^{\a,\b}_{N+2}-\psi^{\a,\b}_N)\chi_N^3\|_{p'}\lesssim N^{\tau(\b,p)} 
 \end{equation*}
where the last inequality can be obtained by a similar argument as in the proof of Lemma \ref{pn}. Again the argument for $D_{N,2}$ is the same as for 
$D_{N,1}$. A similar approach gives $ \|\chi^1_N D_{N,3}\chi^3_N\|_{p\to p} \lesssim N^{\tau(\b,p)-1}\lesssim N^{\tau(\b,p)} $. Hence, we obtain the desired estimates.
\end{proof}

For the remaining estimates of the splitting lemma we need the following result by Hardy-Littlewood and P{\'o}lya \cite[Theorem 3.9]{HLP} that was used in the techniques of Herz \cite{Herz}.
\begin{lemma}[Hardy-Littlewood-P{\'o}lya]\label{he}
Suppose that  $p>1$, and $K$ is an integral operator with kernel $K(x,y)$. Assume that $K(x,y)$ is non-negative, homogeneous of degree $-1$ and  that
\begin{equation}
    \int_0^\infty K(x,1)x^{1/p -1}dx = \int_0^\infty K(1,y)y^{1/{p'}-1}dy <\infty.
\end{equation}
Then, $K$ acts as a bounded operator from $L^p(\R^+,dx)\to L^p(\R^+,dx)$.
\end{lemma}

\begin{theorem}\label{th10.a}
Suppose that $\a,\b\ge -1/2$.
Then there exists $ C > 0 $ independent of $ N $ such that
\begin{equation}\label{hr}
   \|(\chi_{N}^{1}+\chi_{N}^2)D_N (\chi_{N}^{1}+\chi_{N}^2)\|_{p\to p} \le C 
\end{equation}
on $L^p(\R^+,d\mu_\b)$
for all $p$ such that $|\frac{1}{p}-\frac12| < \frac{1}{4(\b +1)}$.
In addition, for all $4/3<p<4$ we also have
\begin{equation}\label{aw}
    \|(\chi_N^{2}+\chi_N^3)D_N(\chi_N^{2}+\chi_N^3)\|_{p\to p} \le C
\end{equation}

\end{theorem}

\begin{proof}
We first discuss \eqref{hr} for the factor $D_{N,1}$.
\begin{align*}
    D_{N,1}f(x)&=\frac{c_{N+1}}{2}   \psi^{\a,\b}_{N+1}(x)\int_0^\infty \frac{f(y)}{(x^2-y^2)}(\psi^{\a,\b}_{N+2}(y)-\psi^{\a,\b}_{N}(y)) y^{2\b+1}dy,
    \\
    &=\frac{c_{N+1}}{2\lambda_N} x^{-\b-1/2}  x^{\b+1/2}\lambda_N^{1/4}\psi^{\a,\b}_{N+1}(x)\int_0^\infty \frac{f(y)}{(x^2-y^2)}y^{\b-1/2}\lambda_N^{3/4}(\psi^{\a,\b}_{N+2}(y)-\psi^{\a,\b}_{N}(y)) y^{\b+3/2}dy.
\end{align*}
Note that $(\chi_N^{1}(x)+\chi_N^2(x))x^{\b+1/2}\lambda_N^{1/4}\psi^{\a,\b}_{N+1}$
and $(\chi_N^{1}(y)+\chi_N^2(y))y^{\b-1/2}\lambda_N^{3/4}(\psi^{\a,\b}_{N+2}(y)-\psi^{\a,\b}_{N}(y))$ are bounded functions. Hence, by changing the  variables $x^2\to x$ and $y^2\to y$ and defining $g(y)=f(\sqrt{y})y^{\frac{\b}{p}}$, we have the following estimates
\begin{align*}
 \int_0^\infty \left| \int_0^\infty x^{-\b/2-1/4}y^{\b/2+1/4-\b/p}  \frac{f(\sqrt{y})y^{\b/p}}{x-y}dy \right|^p x^\b dx=     \int_0^\infty \left| \int_0^\infty \left(\frac{y}{x}\right)^{\frac{\b}{2}-\frac{\b}{p}+\frac{1}{4}} \frac{g(y)}{x-y}dy \right|^p dx
\end{align*}

Following the approach designed by Herz, that is, subtracting the Hilbert transform,  we need to establish the $L^p$ boundedness of the  operator
\begin{align}\label{EqM_m}
   g \rightarrow   \int_0^\infty \left(\left(\frac{y}{x}\right)^{\b/2+1/4-\b/p}-1\right)  \frac{g(y)}{x-y}dy. 
\end{align}
(The above integral, after subtracting the Hilbert transform, is no longer a singular integral operator). Let $$K(x,y)=\left(\left(\frac{y}{x}\right)^{\b/2+1/4-\b/p}-1\right)  \frac{1}{x-y}.$$ Then, by the Hardy-Littlewood-P\'olya inequality, that is Lemma \ref{he}, the $L^p$ boundedness of \eqref{EqM_m}
is valid if 
\begin{equation*}
    \int_0^1 K(1,y)y^{-1/p}dy+\int_1^\infty K(1,y)y^{-1/p}dy<\infty.
\end{equation*}
Note that the first integral is finite whenever 
\begin{equation}\label{Mm_int1}
    \frac\b2+\frac14-\frac\b p-\frac1p>-1, \text{ that is, } \;\frac\b2+\frac54>\frac\b p+\frac1p \text{ or }\;   p>\frac{4\b+4}{2\b+5}, 
\end{equation}
and the second integral is finite whenever 
\begin{equation}\label{Mm_int2}
    \frac\b2+\frac1 4-\frac\b p<\frac1p, \text{ that is, }\;   \frac{4\b+4}{2\b+1} >p.
\end{equation}
Combining \eqref{Mm_int1} and \eqref{Mm_int2} we get
\begin{equation}\label{Mm_int}
   \frac{4\b+4}{2\b+5}  <p<\frac{4\b+4}{2\b+1}.
\end{equation}
To verify bounds for $D_{N,2}$ we consider 
\begin{align*}
   D_{N,2}f(x)&=\frac{c_{N+1}}{2} (\psi^{\a,\b}_{N+2}(x)-\psi^{\a,\b}_{N}(x)) \int_0^\infty \frac{f(y)}{(x^2-y^2)} \psi^{\a,\b}_{N+1}(y)y^{2\b+1}dy,
    \\
    &=\frac{c_{N+1}}{2\lambda_N} x^{-\b+1/2}  x^{\b-1/2}\lambda_N^{3/4}(\psi^{\a,\b}_{N+2}(x)-\psi^{\a,\b}_{N}(x))\int_0^\infty \frac{f(y)}{(x^2-y^2)}y^{\b+1/2}\lambda_N^{1/4}\psi^{\a,\b}_{N+1}(y) y^{\b+3/2}dy.
\end{align*}
Again by changing the variables $x^2\to x$ and $y^2\to y$ and defining $g(y)=f(\sqrt{y})y^{\frac{\b}{p}}$, we reduce the problem to estimating the integral
\begin{align*}
     C \int_0^\infty \left| \int_0^\infty \left(\frac{y}{x}\right)^{\frac{\b}{2}-\frac{\b}{p}-\frac{1}{4}} \frac{g(y)}{x-y}dy \right|^p dx.
\end{align*}
Following an argument similar to the one used to handle $D_{N,1}$ and using Lemma \ref{he}, the $L^p$ boundedness 
is valid if 
\begin{equation}\label{Mm'_Int}
    \frac{4\b+4}{2\b+3}<p<\frac{4\b+4}{2\b-1}
\end{equation}
Thus, $D_{N,1}$ and $D_{N,2}$ together give the following range of $p$ for the uniform boundedness of the partial sums:
\begin{equation}
    \frac{4\b+4}{2\b+3}<p<\frac{4\b+4}{2\b+1}, \text{ that is, }  \left|\frac{1}{p}-\frac12\right|<\frac{1}{4\b+4}.
\end{equation}

Let us now look at the factor $D_{N,3}$ for \eqref{hr}. Note that we can split $D_{N,3}$ into two parts and control each part in a similar way as $D_{N,1}$ and $D_{N,2}$, respectively, since $|\psi^{\a,\b}_{N}(x)|\lesssim N\,\mathcal{H}^{\a,\b}_{N}(x)$. This gives the boundedness of $D_{N,3}$ on  $L^p(d\mu_\b)$ for 
\begin{equation}
    \left|\frac{1}{p}-\frac12\right|<\frac{1}{4\b+4}.
\end{equation}
Thus, combining the estimates for $D_{N,1}, D_{N,2}$ and $D_{N,3}$, we obtain \eqref{hr} for $ \left|\frac{1}{p}-\frac12\right|<\frac{1}{4\b+4}.$
Next we move to an argument for \eqref{aw}. \\

 In dealing with $ (\chi_N^{2}+\chi_N^3)D_{N,1}(\chi_N^{2}+\chi_N^3)$ we write $ D_{N,1}f$ as follows:
\begin{align*}
    & \frac{c_{N+1}}{2}   \psi^{\a,\b}_{N+1}(x)\int_0^\infty \frac{f(y)}{(x^2-y^2)}(\psi^{\a,\b}_{N+2}(y)-\psi^{\a,\b}_{N}(y)) y^{2\b+1}dy,
    \\
    &=\frac{c_{N+1}}{2\lambda_N} \frac{x^{\b}\lambda_N^{1/4}(\lambda_{N}^{1/3}+|x^2-\lambda_{N}|)^{1/4}}{x^{\b}(\lambda_{N}^{1/3}+|x^2-\lambda_{N}|)^{1/4}}\psi^{\a,\b}_{N+1}(x)
    \\ & \quad \times \int_0^\infty \frac{f(y)}{(x^2-y^2)}\frac{y^{\b}\lambda_N^{3/4}(\lambda_N^{1/3}+|y^2-\lambda_N|)^{-1/4}}{(\lambda_N^{1/3}+|y^2-\lambda_N|)^{-1/4}}(\psi^{\a,\b}_{N+2}(y)-\psi^{\a,\b}_{N}(y)) y^{\b+1}dy.
\end{align*}
Note that the following functions are bounded:
$$(\chi_N^{2}(x)+\chi_N^3(x))\, x^{\b}\lambda_N^{1/4} \,(\lambda_N^{1/3}+|x^2-\lambda_N|)^{1/4}\, |\psi^{\a,\b}_{N+1}(x)| \leq C,$$
$$(\chi_N^{2}(y)+\chi_N^3(y))\, y^{\b}\lambda_N^{3/4}\, (\lambda_N^{1/3}+|y^2-\lambda_N|)^{-1/4}\, |(\psi^{\a,\b}_{N+2}(y)-\psi^{\a,\b}_{N}(y))| \leq C. $$ 
Moreover, when $ x, y$ are in the support of  $\chi_N^{2}+\chi_N^3 ,$ they are comparable. Hence, by changing variables  $x^2\to x$ and $y^2\to y$ and taking $g(y)=f(\sqrt{y})y^{\frac{\b}{p}}$, we reduce the problem to 
estimating the following integral:
\begin{align}\label{X_2+X_3M_m}
     \int_{\sqrt{N/2}}^{\sqrt{2N}} \Big| \int_{\sqrt{N/2}}^{\sqrt{2N}} \Big(\frac{\lambda_N^{1/3}+|y-\lambda_N|}{\lambda_N^{1/3}+|x-\lambda_N|}\Big)^{\frac{1}{4}} \frac{g(y)}{x-y}dy \Big|^p dx.
\end{align}
Again, by changing the variables $\lambda_N-x=X$ and $\lambda_N-y=Y$ and noticing that
\begin{equation*}
    \Big(\frac{\lambda_N^{1/3}+|Y|}{\lambda_N^{1/3}+|X|}\Big)^{\frac{1}{4}}\lesssim \Big(\frac{\lambda_N^{1/3}}{\lambda_N^{1/3}+|X|}\Big)^{\frac{1}{4}}+\Big(\frac{|Y|}{\lambda_N^{1/3}+|X|}\Big)^{\frac{1}{4}}\lesssim \Big(\frac{\lambda_N^{1/3}}{\lambda_N^{1/3}+|X|}\Big)^{\frac{1}{4}}+\Big(\frac{|Y|}{|X|}\Big)^{\frac{1}{4}},
\end{equation*}
 the problem is reduced to estimating the integrals
     $$ \int_{\lambda_N -\sqrt{N/2}}^{\lambda_N -\sqrt{2N}} \Big| \int_{\lambda_N -\sqrt{N/2}}^{\lambda_N -\sqrt{2N}}  \left(\big(\frac{|Y|}{|X|}\big)^{\frac{1}{4}}+\big(\frac{\lambda_N^{1/3}}{\lambda_N^{1/3}+|X|}\big)^{\frac{1}{4}} \right)\frac{\tilde{g}(Y)}{X-Y}dY \Big|^p dX .$$
The boundedness of the Hilbert transform along with the idea of Herz shows that the above is bounded by a constant times $\| f\|_{L^p(d\mu_\b)}.$

A similar argument for $D_{N,2}$ with $\chi_N^2+\chi_N^3$ gives the boundedness on $L^p(d\mu_\b)$ for $p>\frac43$, and hence $D_{N,1}$ and $D_{N,2}$ together give the boundedness for $\frac{4}{3}<p<4.$
The argument for $D_{N,3}$ also holds similarly, giving boundedness on $L^p(d\mu_\b)$ for $\frac43<p<4$.
Again, combining the boundedness  of $D_{N,1}, D_{N,2}$ and $D_{N,3}$, we obtain \eqref{aw} for $\frac{4}{3}< p<4$.

\end{proof}

\section{Sufficient condition}\label{SME}
The following theorem is  the core of our project. By providing a sufficient condition for the uniform boundedness of the Bochner-Riesz means, it forms the  analogue of Theorem~\ref{verlp}, where we established  a necessary condition.  The proof of this sufficient condition is considerably more intricate and requires an elaborate argument. 

\begin{theorem}\label{Th:Sufficient}
For any   $1 < p < \infty$ suppose that 
$$
\delta >  \max\left\{0,\quad \frac{2}{3}\left|\frac{1}{2}-\frac{1}{p}\right|-\frac{1}{6}, \quad
2(\b+1)\left|\frac{1}{2}-\frac{1}{p}\right|
-\frac{1}{2}
\right\}
.$$
Then  the  Bochner-Riesz means of order $\delta$ for the operator $\Delta_{\a,\b}$ are uniformly bounded on $L^p,$ i.e.
$$
\sup_{N>0}\| S_N^{\delta}(\Delta_{\a,\b})\|_{p \to p} \le C < \infty.
$$

\end{theorem}

The proof of this theorem is organised into several lemmas and propositions, and draws on ideas developed in the preceding sections. Before presenting the technical details, we outline the overall strategy and structure of the argument. We employ the decomposition from subsection \ref{Subsec:SplitLem}, as formulated in Lemma \ref{elem1},  $T=S_N^{\delta}(\Delta_{\a,\b})$. The component associated with $\chi^4_N$ is elementary and follows directly from Lemma \ref{2N}.

The component $(\chi_N^{1}+\chi_N^{2})T(\chi_N^{1}+\chi_N^{2})$ corresponds to the Bessel operators and is linked to \cite[Part III, Section 4]{COSY}; see also \cite{Herz}. We show that this contribution is bounded above the region $B$ defined in \eqref{b}.   The spectral multiplier estimates on the Bessel-type part--$(\chi^1_N+\chi^2_N)T(\chi^1_N+\chi^2_N)$ are proved in Theorem~\ref{th5.4} below which follows from Lemma~\ref{l5.4} and Proposition \ref{PropFxi}. Moreover,  the corresponding result on the uniform boundedness of the Bochner-Riesz means is given in Remark~\ref{Rmk:chi1chi2} below.

The term $(\chi_N^{2}+\chi_N^3)T(\chi_N^{2}+\chi_N^3)$ is called the harmonic oscillator type part and it extends the estimates above region $C$ defined in \eqref{c}. Its proof relies on a transplantation phenomenon, formulated in Proposition \ref{trans} below, together with the uniform boundedness of the Bochner–Riesz means for the one-dimensional Hermite expansions established in \cite{Th89}.

For $p\le 2$, the term $\chi_N^{1}T\chi_N^3$ can be absorbed into the contribution $(\chi_N^{2}+\chi_N^3)T(\chi_N^{2}+\chi_N^3)$, while $\chi_N^{3}T\chi^1_{N}$ aligns more closely with the approach used for $(\chi_N^{1}+\chi_N^{2})T(\chi_N^{1}+\chi_N^{2})$, see Corollary~\ref{TT}. The factor $\chi^2_N$ does not affect the arguments employed for either $(\chi_N^{2}+\chi_N^3)T(\chi_N^{2}+\chi_N^3)$ or $(\chi_N^{1}+\chi_N^{2})T(\chi_N^{1}+\chi_N^{2})$. Its role is to separate $\chi^1_N$ from $\chi^3_N$, and this separation is crucial for the argument used in Lemma \ref{lemchi3chi1} below.

Note that our approach to the component $(\chi_N^{1}+\chi_N^{2})T(\chi_N^{1}+\chi_N^{2})$ is based on the finite speed of propagation property for the corresponding wave equation.  Hence, it is necessary to replace $\Delta_{\a,\b}$ by $\sqrt{\Delta_{\a,\b}}$ and hence $ N $ needs to be replaced by 
 $N^2$. In fact, for this component and $\chi^3_N T\chi^1_N$, instead of directly getting the boundedness of the Bochner-Riesz means, we obtain bounds on their spectral multipliers for $T=F(\sqrt{\Delta_{\a,\b})}$. In these cases, the latter approach is easier to follow, and it does imply the boundedness of the Bochner-Riesz means for these components. The result for the latter part is mentioned in Remark~\ref{Rmk:chi3chi1}\\

With the next lemma, we begin  the analysis of the component $(\chi_N^{2}+\chi_N^3)T(\chi_N^{2}+\chi_N^3).$ Recall that 
\begin{equation*}
    M_{\gamma} f(x) = x^{\gamma} f(x).
\end{equation*}

\begin{lemma}\label{elem2}
   For $\sigma>0$  let $d\omega_\sigma=x^\sigma dx$, and let $T$ be a linear operator acting on $L^p(\R^+)$. Then for any $0<a<b$ such that $|b/a|\le K$, for some fixed $K>0$, and for all  $1 \le p \le \infty$
   \begin{align*}
   K^{-\sigma/p}\big\|\chi_{[a,b]}T\chi_{[a,b]}\big\|_{L^p \to L^p} \le  \big\|\chi_{[a,b]}T\chi_{[a,b]}\big\|_{L^p(d\omega_\sigma) \to L^p(d\omega_\sigma)} \le K^{\sigma/p}\big\|\chi_{[a,b]}T\chi_{[a,b]}\big\|_{L^p \to L^p}.
\end{align*}  
Moreover, for any positive measure $ \mu$ on $\R^+$ and $\gamma\in\RR,$
\begin{align*}
   K^{-|\gamma|}\big\|\chi_{[a,b]}T\chi_{[a,b]}\big\|_{L^p(d\mu) \to L^p(d\mu)} \le  \big\|\chi_{[a,b]}M_\gamma T M_\gamma^{-1}  \chi_{[a,b]}\big\|_{L^p(d\mu) \to L^p(d\mu)} \le K^{|\gamma|}\big\|\chi_{[a,b]}T\chi_{[a,b]}\big\|_{L^p(d\mu) \to L^p(d\mu)}
\end{align*}  
for all $1 \le p \le \infty$.
\end{lemma}

\begin{proof} Defining $ f_{a,b}= f \chi_{[a,b]} $ and $ T_{a,b} = \chi_{[a,b]} T\chi_{[a,b]} ,$ we have
$$ \int_0^\infty |\chi_{[a,b]} T\chi_{[a,b]}f(x)|^p x^\sigma\, dx  \leq  (b/a)^\sigma  \| \chi_{[a,b]}T \chi_{[a,b]}\|_{L^p \rightarrow L^p} \int_0^\infty |f(x)|^p x^\sigma dx.$$ 
This proves half of the first inequality. The other half is proved similarly.  Note that the first set of inequalities is true even if $ L^p(dx) $ is replaced by $ L^p(d\mu) $ for any positive measure $ \mu.$ For the second set of inequalities use the easy to check estimates 
$
\|\chi_{[a,b]} M_\gamma\|_{L^p(d\mu) \to L^p(d\mu)} \le b^{\gamma} $ valid for $ \gamma >0 $ and 
$ \|\chi_{[a,b]} M_\gamma\|_{L^p(d\mu) \to L^p(d\mu)} \leq  a^{-\gamma}$
valid for $ \gamma <0.$ Hence,  the second set of inequalities follows  immediately from the first set. 

\end{proof}

 For later use we record the following corollary and a transplantation theorem.

\begin{coro}\label{co1}
    Let $\sigma>0$, $1 \le p \le 2$, $d\omega_\sigma=x^\sigma dx$,  $T$ is a linear operator acting on $L^p(\R^+)$. Then 
    $$
    \big\|\chi_{[0,b]}M_{\sigma/2}^{-1}TM_{\sigma/2}\chi_{[a,b]}\big\|_{L^p(d\omega_\sigma) \to L^p(d\omega_\sigma)}\le K^{\sigma(1/p-1/2)}\big\|T\big\|_{L^p \to L^p}.
    $$
    \end{coro}
\vskip0.15truein

\begin{proposition}[Transplantation] \label{trans}
Suppose that $\a_1,\a_2\ge -1/2 $. 
Then for any bounded Borel function  $F$ on $\R^+$
\begin{equation}\label{pos1}
    \big\|F({\tilde{\Delta}_{\a_1}}) \big\|_{p\to p} \sim 
\big\|F({\tilde{\Delta}_{\a_2}})\big\|_{p\to p}
\end{equation}
for all $1< p <\infty$.
\end{proposition}
\begin{proof} The relation \eqref{pos1} is essentially equivalent to \cite[Theorem 7.1.2]{ThBook}. Indeed, if 
\begin{equation*}
    \tau^{\a_1}_{\a_2}f=\sum_{k=0}^\infty\langle f,h_k^{\a_2}\rangle h_k^{\a_1},
\end{equation*}
 is the transplantation operator associated to  the Laguerre functions of Hermite type (see \eqref{eqDelta2}), then 
    $\tau^{\a_1}_{\a_2}=(\tau^{\a_2}_{\a_1})^{-1}$ is bounded on 
    all $L^p$, see \cite[Theorem 7.1.2]{ThBook}.
\end{proof}
This transplantation theorem leads us to the following theorem which gives the sufficient condition for the uniform boundedness  of the Bochner-Riesz means associated to the operator $\tilde{\Delta}_{\a}$.

\begin{theorem}\label{verlp2po}
Let $ \a \geq -1/2.$ For $1< p< \infty,$ let
 $$
\delta >  \max \left\{0,\quad \frac{2}{3}\left|\frac{1}{2}-\frac{1}{p}\right|-\frac{1}{6} 
\right\}.
$$
 Then the  Bochner-Riesz means of order $\delta$ for  $\tilde{\Delta}_{\alpha}$ are  uniformly bounded on $L^p$, that is
$$
\sup_{R>0}\| S_R^{\delta}(\tilde{\Delta}_{\alpha})\|_{p \to p} \le C < \infty.
$$
\end{theorem}
\begin{proof}
By choosing $F(\lambda) = ( 1-\lambda/R)_+^\delta$, $\a_1=\a$ and $\a_2=\pm\frac{1}{2}$ in Proposition~\ref{trans}, we obtain 
   \begin{equation}
    \big\|S_R^{\delta}({\tilde{\Delta}_{\a}}) \big\|_{p\to p} \sim 
\big\|S_R^{\delta}({\tilde{\Delta}_{\pm 1/2}})\big\|_{p\to p}
\end{equation}
for all $1< p <\infty$. Now, by \cite[Theorem 2.1]{Th89} and a transplantation argument, $S_R^{\delta}({\tilde{\Delta}_{\pm 1/2}})$ is uniformly bounded on $L^p$ for $\delta>\frac{2}{3p}-\frac{1}{2}$ , if $p<2$. See also \cite[Theorem 6.4.1]{ThBook}. Note that $\a=\pm 1/2$ corresponds to the action of the operator $-\frac{d^2}{dx^2}+x^2$ on even or odd functions on the real line. This gives us the desired result.

\end{proof}

\begin{theorem} \label{TransTh}
Let $F$ be a bounded Borel function on $\R^+$. For any $0<a<b$ such that $|b/a|\le K$, for some fixed $K>0$, and for all  $1 < p < \infty$ we have
 \begin{eqnarray*}
     \big\|\chi_{[a,b]}F(\Delta_{\a,\b})\chi_{[a,b]}\big\|_{L^p(\mu_\beta) \to L^p(\mu_\beta)}\le 
          \big\|F({\tilde{\Delta}_{\a}}) \big\|_{L^p(\R^+) \to L^p(\R^+)} \sim  C\big\|F(-dx^2+x^2)\big\|_{L^p(\R^+) \to L^p(\R^+)}.
 \end{eqnarray*}  
 Additionally if $p\le 2$ then we also have
  \begin{eqnarray*}
     \big\|\chi_{[0,b]}F(\Delta_{\a,\b})\chi_{[a,b]}\big\|_{L^p(\mu_\beta) \to L^p(\mu_\beta)}\le 
          \big\|F({\tilde{\Delta}_{\a}}) \big\|_{L^p(\R^+) \to L^p(\R^+)} \sim  C\big\|F(-dx^2+x^2)\big\|_{L^p(\R^+) \to L^p(\R^+)}.
 \end{eqnarray*} 
\end{theorem}

\begin{proof}
Note that by \eqref{unitary1}
\begin{equation}\label{ab1}
F(\Delta_{\a,\beta})f(x)= M^{-1}_{1/2+\beta}F(\tilde{\Delta}_{\alpha})M_{1/2+\beta} f(x). 
\end{equation}
Hence, by Lemma \ref{elem2} 
$$
\big\|\chi_{[a,b]} F(\Delta_{\a,\beta})\chi_{[a,b]}
\big\|_{L^p(\mu_\beta) \to L^p(\mu_\beta)}
\sim \big\| \chi_{[a,b]}F(\tilde{\Delta}_{\alpha})\chi_{[a,b]}\big\|_{L^p(\R^+) \to L^p(\R^+)}
$$
so that
$$
\big\|\chi_{[a,b]} F(\Delta_{\a,\beta})\chi_{[a,b]}
\big\|_{L^p(\mu_\beta) \to L^p(\mu_\beta)}
\leq C \big\| F(\tilde{\Delta}_{\alpha})\big\|_{L^p(\R^+) \to L^p(\R^+)}.
$$
Moreover, for $p\le 2$, by Corollary \ref{co1}
$$
\big\|\chi_{[0,b]} F(\Delta_{\a,\beta})\chi_{[a,b]}
\big\|_{L^p(\mu_\beta) \to L^p(\mu_\beta)}
\leq C \big\| F(\tilde{\Delta}_{\alpha})\big\|_{L^p(\R^+) \to L^p(\R^+)}.
$$
Next, by the transplantation result in Proposition~\ref{trans},
$$
\big\| F(\tilde{\Delta}_{\alpha})\big\|_{L^p(\R^+) \to L^p(\R^+)}
\sim 
\big\| F(\tilde{\Delta}_{-1/2})\big\|_{L^p(\R^+) \to L^p(\R^+)}=\big\|F(-dx^2+x^2)\big\|_{L^p(\R^+) \to L^p(\R^+)}.
$$
This completes the proof. 
\end{proof}

 By specialising to the function $ F(\lambda) = (1-\lambda/R)_+^\delta $ we obtain the following results for the Bochner-Riesz means:

\begin{coro}\label{Trans}
     For any $0<a<b$ such that $|b/a|\le K$, for some fixed $K>0$, and for all  $1 < p < \infty$
$$
\sup_{R>0}\big\|\chi_{[a,b]} S_R^{\delta}(\Delta_{\a,\beta})\chi_{[a,b]}\big\|_{L^p(d\mu_\b) \to L^p(d\mu_\b)}
\le C\sup_{R >0}\big\|S_R^{\delta}(-dx^2+x^2)\big\|_{L^p(\R^+) \to L^p(\R^+)}.
 $$  
Moreover for $ 1 < p \le 2 $
 $$
\sup_{R >0}\big\|\chi_{[0,b]} S_R^{\delta}(\Delta_{\a,\beta})\chi_{[a,b]}\big\|_{L^p(d\mu_\b) \to L^p(d\mu_\b)}
\le C\sup_{R >0}\big\| S_R^{\delta}(-dx^2+x^2)\big\|_{L^p(\R^+) \to L^p(\R^+)}.
$$
\end{coro}

\begin{coro}\label{TT}
Assume that  $0<a<b$ are such that $|b/a|\le K$, for some fixed $K>0.$ Then for any $ 1 < p <\infty,$ we have
$$
\sup_{R >0}\big\|\chi_{[a,b]} S_R^{\delta}(\Delta_{\a,\beta})\chi_{[a,b]}\big\|_{L^p(d\mu_\b) \to L^p(d\mu_\b)}
\le C
 $$ 
 provided $ \delta $ satisfies
 $$
\delta >  \max\left\{0,\quad \frac{2}{3}\left|\frac{1}{2}-\frac{1}{p}\right|-\frac{1}{6} 
\right\}.
$$
If in addition  $ 1 < p\le 2,$ we also have
\begin{equation}\label{chi1to3}  
\sup_{R >0}\big\|\chi_{[0,b]}S_R^{\delta}(\Delta_{\a,\beta})\chi_{[a,b]}\big\|_{L^p(d\mu_\b) \to L^p(d\mu_\b)}
\le C
\end{equation}
for the same range of $\delta$ as above. In particular, for $1 < p\le 2 $ and $ N \in \mathbb N,$ we have
\begin{equation}\label{3chi} 
\sup_{R >0}\big\|\chi^1_N S_R^{\delta}(\Delta_{\a,\beta})\chi_N^3\big\|_{L^p(d\mu_\b) \to L^p(d\mu_\b)}
\le C
\end{equation}
for the same range of $\delta$ with a constant independent of $ N.$
\end{coro}

\begin{proof} 
By \cite[Theorem 5.5]{Th89} and Corollary~\ref{Trans} one can obtain the first two estimates. Now, taking $b=2N$ and $a=N/4$ in \eqref{chi1to3}, we get
    $$
\sup_{R >0}\big\|(\chi^1_N+\chi_N^2+\chi_N^3) S_R^{\delta}(\Delta_{\a,\beta})(\chi^2_N+\chi_N^3)\big\|_{L^p(d\mu_\b) \to L^p(d\mu_\b)}
\le C.
$$
Clearly, 
$$
\big\|\chi^1_N S_R^{\delta}(\Delta_{\a,\beta})\chi_N^3\big\|_{p\to p}\le 
\big\|(\chi^1_N+\chi_N^2+\chi_N^3) S_R^{\delta}(\Delta_{\a,\beta})(\chi^2_N+\chi_N^3)\big\|_{p\to p}
$$
which ends the proof of the above estimate.
\end{proof}

\subsection{Bessel-type part and its interaction with the harmonic oscillator part}\label{sec5.1}

In the part corresponding to the component
 $$(\chi_N^1+\chi_N^2)F(\Delta_{\alpha,\beta})(\chi_N^1+\chi_N^2),$$ the spectral multipliers behave like Bessel-type operators.
This is because in this region the harmonic oscillator
potential \(x^2\) is not crucial at the relevant spectral scale, while the singular part
\[
-\frac{d^2}{dx^2}
-\frac{2\beta+1}{x}\frac{d}{dx}
+\frac{\alpha^2-\beta^2}{x^2}
\]
is exactly the Bessel-type inverse-square operator considered in Part III, Section 4 of
\cite{COSY}. With \(n=2\beta+2\) and \(c=\alpha^2-\beta^2\), the restriction estimate for
\(\Delta_n+c/r^2\) gives the local \(L^p\to L^2\) estimate which yields the contribution
\[
\delta>2(\beta+1)\left|\frac1p-\frac12\right|-\frac12.
\]
This is precisely the second term in the sufficient condition. This observation can be expressed using the notion of spectral multipliers.\\ 

We plan to use some properties of the wave operator in studying  $(\chi_N^1+\chi_N^2)F_\delta(\Delta_{\alpha,\beta})(\chi_N^1+\chi_N^2),$  where $ F_\delta(\lambda) = (1-\lambda/N)_+^\delta $ so that  $ S_N^\delta(\Delta_{\a,\b}) = F_\delta(\Delta_{\a,\b}). $ This will be handled by proving a general multiplier theorem for the operators $ F(\sqrt{\Delta_{\a,\b}}) .$ For technical reasons, in studying such multipliers  we will be assuming that $ F $ is supported on the interval $[N/2,N] .$  This necessitates the splitting of the multiplier $ F_\delta(\lambda) = F_\delta^{(1)}(\lambda)+F_\delta^{(2)}(\lambda) $  with $ F_\delta^{(1)}(\lambda) = \eta(\lambda/N)F_\delta(\lambda)$ where $ \eta $ is a smooth function supported in a small neighborhood of  $[0,1/4] $ and $ \eta(\lambda) =1 $ on $[0,1/4].$ The part of the Bochner-Riesz means corresponding to $ F_\delta^{(1)}(\Delta_{\a,\b})$ can be taken care of by known techniques since the multiplier is smooth and compactly supported, see the proposition below.  The remaining part  is
 $$ F_\delta^{(2)}(\Delta_{\a,\b}) = ( 1-\eta(\Delta_{\a,\b}/N))\, (1-\Delta_{\a,\b}/N)_+^\delta  $$ 
 and note that $ F_\delta^{(2)}$ is supported on $[N/4,N].$ Replacing $ N $ by $ N^2 $ we are interested in the uniform boundedness of 
$$(\chi_{N^2}^1+\chi_{N^2}^2)F_\delta^{(2)}(\Delta_{\alpha,\beta})(\chi_{N^2}^1+\chi_{N^2}^2)$$ 
or equivalently, in the uniform boundedness of the operators 
$$(\chi_{N^2}^1+\chi_{N^2}^2)F_\delta^{(3)}(\sqrt{\Delta_{\alpha,\beta}})(\chi_{N^2}^1+\chi_{N^2}^2)$$
where $ F_\delta^{(3)}(\lambda) = F_\delta^{(2)}(\lambda^2) $ which is supported on $ [N/2,N].$ These operators are handled in Theorem \ref{th5.4}.\\

Let us take care of $(\chi_N^1+\chi_N^2)F_\delta^{(1)}(\Delta_{\alpha,\beta})(\chi_N^1+\chi_N^2) $ before we  move on to Theorem \ref{th5.4}. Let $  m(\lambda) $ be a compactly supported smooth function and define $ m_N(\lambda) = m(\lambda/N).$ Then it is enough to prove the following proposition.

 \begin{proposition}\label{Prop:mN} The operators $ m_N(\Delta_{\a,\b}) $ are uniformly bounded on $ L^p(\R^+, d\mu_\b) $ for all $ 1 \leq p \leq \infty.$
\end{proposition}

We deduce the above proposition from  the following result proved in \cite{DOS}. In this work the authors have considered multipliers associated to self adjoint, positive definite operators $ A $ of order 2 acting on $ L^2(X, d\mu) $ where $ \mu $ is assumed to have the doubling property.  The operator $-A $ generates a  semi group $ e^{-t A} $ whose  kernel  is denoted by $ K_{e^{-tA}}(x,y) .$ They assume that this heat kernel satisfies the Gaussian upper bounds 
\begin{equation}\label{gaussian}
K_{e^{-tA}}(x,y) \leq C\, \mu(B(y, \sqrt{t}))^{-1}\, e^{-\frac{b}{t} \rho(x,y)^2},
\end{equation}
where $\rho(x,y)$ denotes the distance on the metric space $X$, and $C, b$ are some positive constants.\\
The following result has been proved in \cite[Theorem 6.2]{DOS}.

\begin{theorem}\label{Thm:DOS}
	Suppose that $A$ is as above whose heat kernel satisfies the condition \eqref{gaussian}. Then any Borel function $F$ such that $\supp F\subset[0,N]$ the kernel $ K_{F(A)} $ of the operator $ F(A) $ satisfies 
	\begin{equation}\label{eq:DOS}		\|K_{F(A)}(\cdot,y)\|^2_{L^2(X, d\mu)}\le C \mu(B(y, N^{-1/2}))^{-1}\| \delta_NF\|_{\infty}^2.
	\end{equation}
	 For any $s>d/2, d$  being the dimension of $ X $, we also have
	\begin{equation}\label{eq2:DOS}
		\|F(A)\|_{L^q(X, d \mu)\to L^q(X,d \mu)}\le C\|\delta_NF\|_{W^{s,\infty}}
	\end{equation}
	for all  $q\in[1,\infty]$.
\end{theorem}
{\bf Proof of Proposition \ref{Prop:mN}:} With $X=\R^+, \mu= \mu_\beta$  and $A=\Delta_{\a,\b} $ we apply  Theorem~\ref{Thm:DOS} to the operator $ F(\Delta_{\a,\b}) $ where  $F(\lambda)=m_N(\lambda)=m(\lambda/N)$ with $\supp m\subset [0,1/2].$ We only need to verify that the heat kernel associated to  $\Delta_{\a,\b}$ satisfies the Gaussian upper bound \eqref{gaussian}. But this is easy to see. Indeed, the functions $ \psi_k^{\a,\b} $ are eigenfunctions of $ \Delta_{\a,\b} $ with eigenvalues $ \lambda_k = (4k+2\alpha+2)$ the heat kernel is explicitly given by 
$$ K_{e^{-tA}}(x,y) = 2(xy)^{\a-\b} \, \sum_{k=0}^\infty e^{-(4k+2\a+2)t}\, \frac{\Gamma(k+1)}{\Gamma(k+\a+1)}\, L_k^\a(x^2)\, L_k^\a(y^2)\, e^{-\frac{1}{2}(x^2+y^2)}.$$
Using the generating function identity (1.1.47) in \cite{ThBook} and estimating the Bessel function we can get the Gaussian bound \eqref{gaussian}.
Since $m$ is a smooth and compactly supported function, we have $\|m\|_{W^{s,\infty}}<\infty$ for all $s\ge 0$. Therefore, from \eqref{eq2:DOS}, we get
	\begin{equation*}
	\|m_N(\Delta_{\a,\b})\|_{L^q(\R^+,d\mu_\b)\to L^q(\R^+,d\mu_\b)}\le C\|m\|_{W^{s,\infty}}<\infty
\end{equation*}
uniformly in $N$ for all $q\in[1,\infty]$.\\



We will now state Theorem \ref{th5.4} and work towards a proof. Let us recall that for $ s > 0,$ the Sobolev space $ H^s= H^s(\R)$ is the subspace of $ L^2(\R)$ consisting of those $ \varphi $ for which the following norms are finite:
$$ \| \varphi \|_{H^s}^2 = \int_{-\infty}^\infty |\widehat{\varphi}(y)|^2\, (1+y^2)^{s}\, dy.$$


\begin{theorem}\label{th5.4}
For  $\a,\b \ge  -1/2$  and 
 $1 \le p \le 2$ define
$$
s_p(\a,\b) \coloneqq  \max\left\{\frac12, \quad
2(\b+1)\left|\frac{1}{2}-\frac{1}{p}\right|
\right\}.
$$
For any  bounded Borel
function 
$F: [0, \infty) \to \mathbb{C}$ with  $\supp\, F \subset [N/2,N] $ define $ \delta_NF(t) = F(Nt).$
Then  for any $ s > s_p(\a,\b) $ there exists $ C $ independent of $ N $ such that 
\begin{eqnarray*}
	\|(\chi^1_{N^2} +\chi^2_{N^2})  F\Big({\sqrt{\Delta_{\a,\b}}}\Big) (\chi^1_{N^2} +\chi^2_{N^2}) \|_{L^r(d\mu_\b)\to L^r(d\mu_\b)}\leq    C\, \| \delta_N F\|_{H^{s}} 
\end{eqnarray*}
for all $ p \le r \le p^\prime.$ Thus if  $ \delta_NF \in H^{s}, s > s_p(\a,\b),$ with $ \| \delta_N F\|_{H^{s}} $ uniformly bounded, then the operators 
$$ (\chi^1_{N^2} +\chi^2_{N^2})  F\big(\sqrt{\Delta_{\a,\b}}\big) (\chi^1_{N^2} +\chi^2_{N^2}) $$ are uniformly bounded on $ L^p.$
\end{theorem}

\begin{remark}\label{Rmk:BRmeans}
	 In most of the applications, e.g. Bochner-Riesz means, we have $ F(\lambda) = G(\lambda/N) $ so that $ \delta_NF =G $ and hence $ \| \delta_NF \|_{H^s} $ is independent of $N.$ In such a situation we only need to check if $ G = \delta_NF \in H^s $ or not.
\end{remark}

\begin{remark}\label{Rmk:chi1chi2}
In Theorem~\ref{th5.4}, if we choose $F(\lambda)=(1-\lambda^2/N^2)^\delta_+$, then $F\in H^s$ if and only if $\delta>s-1/2$, and from Remark~\ref{Rmk:BRmeans}, $ \| \delta_NF \|_{H^s} $ is independent of $N$.  Therefore, for $S^\delta_{N^2}(\Delta_{\a,\b})=F(\sqrt{\Delta_{\a,\b}})$, we obtain the uniform boundedness of the Bochner-Riesz means corresponding to the Bessel-type part of the operator. Indeed, for $\a,\b\ge -1/2$, $1\le p\le \infty$, and 
$$
\delta >  \max\left\{0, \quad
2(\b+1)\left|\frac{1}{2}-\frac{1}{p}\right|
-\frac{1}{2}
\right\}
$$
the part of Bochner-Riesz means corresponding to $(\chi^1_{N^2}+\chi^2_{N^2})$ of order $\delta$  is uniformly bounded on $L^p
(\R^+,x^{2\b+1}dx)$, that is,
$$
\sup_{N>0}\|(\chi^1_{N^2} +\chi^2_{N^2})S_{N^2}^{\delta}( \Delta_{\a,\b})(\chi^1_{N^2} +\chi^2_{N^2})\|_{p \to p} \le C < \infty.
$$
\end{remark}

 In order to prove the above theorem we require  Lemma \ref{l5.4} and Proposition \ref{PropFxi} and we now proceed to prove them. In proving the lemmas, we split the operator $F(\sqrt{\Delta_{\a,\b}})$ into two parts as follows. 
$$ F\Big(\sqrt{\Delta_{\a,\b}}\Big)=(F*\xi)\Big(\sqrt{\Delta_{\a,\b}}\Big)+(F-F*\xi)\Big(\sqrt{\Delta_{\a,\b}}\Big),$$
where $\xi$ is a smooth, compactly supported function on $\RR$ satisfying the assumptions stated in Proposition~\ref{wuj1a} below. The part $(F-F*\xi)(\sqrt{\Delta_{\a,\b}})$ is considered in Lemma~\ref{l5.4} and the remaining part $ (F*\xi)(\sqrt{\Delta_{\a,\b}})$ in Proposition~\ref{PropFxi}. 
The operator $F* \xi$ separates the smooth spectral part from the discrete spectral irregularities that come from the $x^2$ potential. The part $(F-F*\xi)$, which corresponds to discrete spectral irregularities, is not essential; nevertheless, it is not negligible. 
 In order to estimate these operators, we consider the following norm, introduced in \cite{CowS, DOS}, which is essential in situations where the operators have discrete spectrum.
For an even Borel function $F$
with  $\supp F\subseteq [-1, 1]$ we define
the norm $\|F\|_{N,q}$ by
\begin{equation}\label{Def_N,2}
\|F\|_{N,q}\coloneqq\Bigg({1\over 2N}\sum_{\ell=1-N}^{N} \sup_{\lambda\in
 [{\ell-1\over N}, {\ell\over N})} |F(\lambda)|^q\Bigg)^{1/q},  
\end{equation}
 where $q\in [1, \infty)$ and $N\in \NN$. For $q=\infty$, we put
 $\|F\|_{N, \infty}=\|F\|_{{\infty}}$.
 It is obvious that $\|F\|_{N,q}$ increases monotonically in $q.$ For our purposes in this work we only need the $q=2$ version. We also need the following proposition from \cite{DOS} for $p=2$.\\
 
 \begin{proposition}\label{dos}
 Suppose that $\xi \in C_c^{\infty}$ is a function such that
{\rm supp}~$\xi \subset [-1,1]$,  $\xi \ge 0$,   $\widehat{\xi}(0)=1$   and
$\widehat{\xi}^{(k)}(0)=0$  for all $1\le k  \le [s]+2$. Set $\xi_N(\lambda)=N\xi(N\lambda).$  Then for any Borel function $ G $ such that
{\rm supp} $G\subset [0,1]$  and $ N \geq 2,$ we have
$$\|G-G*\xi_N\|_{N,2} \le CN^{-s}\|G\|_{H^s}$$
for all $s>1/2.$
 \end{proposition}
 \begin{proof} We refer to  \cite[Proposition 4.6]{DOS}. Here is a sketch of the proof.
 We first observe that the support of $ G\ast\xi_N $ is contained in $ [-1/N, 1+1/N]\subset [-1,2]$ for $ N \geq1$. Further,  for any $ \lambda \in [-1, 2]$ we have
 $$ |\xi_N\ast G( \lambda)|^2 \leq  \|\xi_N\|_2^2\, \int_{\lambda-1/N}^{\lambda+1/N} |G(t)|^2\ dt.$$  
 Consequently, we have
 \begin{align}\label{f: Np estimate  2}
 	\| \xi_N*G\|_{N,2}&=\Bigg(\frac{1}{N}\sum_{i=1-N}^{N}\sup_{\lambda\in[\frac{i-1}{N},\frac{i}{N})}|\xi_N* G(\lambda)|^2\Bigg)^{1/2}\notag\\
 	&\le\frac{\|\xi_N \|_{2}}{N^{1/2}}\Bigg(  N\, \sum_{i=1-N}^{N}\int_{-1/N+(i-1)/N}^{1/N+i/N}|G(t)|^2 dt\Bigg)^{1/2}\notag\\
 	&\le  \|\xi_N \|_{2}\, \|G\|_{2}.
 \end{align} 

 Now define the function $ \zeta_s$ by the relation
$
\widehat{\zeta_s}(\lambda) = \big(1-\widehat\xi(\lambda)\,\big) |\lambda|^{-s}.$
Observe first that
\begin{equation}\label{f: subclaim}
  \sum_{i\in\Z} \sup_{t\in[i-1,i]}  |\zeta_s  *  G    |^{2} \leq C   \|
         G\|^2_{2}, \qquad\forall G \in L^2(\R) .
\end{equation}
Indeed, Fourier analysis  shows that  $| \zeta_s(t)|\leq C_1  |t|^{s-1}$
when $|t|  \leq 1$ and  $| \zeta_s(t)|  \leq C_2  |t|^{s-k-1}$ when $|t|
\geq   1$.  Therefore,   we   may  write  $\zeta_s$  as    $\sum_{j\in\Z}
\xi_{s,j}({\cdot}-j)$,  where  $\supp \xi_{s,j}  \subseteq [-1,1]$ and
$\sum_{j\in\Z} \|  \xi_{j,s} \|_{2} <  \infty$ (this is  where we
require that $s > 1/2$). The argument leading to  \eqref{f: Np estimate 2} proves \eqref{f: subclaim}.
The proof of our claim is now straightforward. 
Indeed, writing $ G = \delta_NF $ so that $ G -G \ast\xi_N= \delta_N(F- F \ast \xi)$ we have
\begin{align*}
\|\delta_N (F-\xi \ast F)\|_{N,2}^2  &   =     {N^{-1}}
\sum_{i=1-N}^{N} \sup_{t\in  [\frac{i-1}{N},\frac{i}{N}]} |  (F  -
\xi* F)(Nt) |^{2}  \\[4pt] & \leq  N^{-1} \sum_{i=-\infty}^{\infty}
\sup_{t\in[i-1,i]} | \zeta_s * I_s F(t) |^2,
\end{align*}
where  $\zeta_s$  is as   above   and  $\widehat{(I_s F)}  =   |\cdot|^{s}
\widehat{F}$.  Therefore, by appealing to \eqref{f: subclaim} and making use of the relation between $ \| F\|_{H^s} $ and $ \| \delta_NF\|_{H^s}$ we get
\begin{eqnarray}\label{f: EstF-Fxi}
\|\delta_{N}(F - \xi* F)\|_{N,2}  \leq C \, N^{-1/2}
\|  I_s F\|_{2}  \leq   C \, N^{-s} \|\delta_{N}F\|_{H^s}.
\end{eqnarray}
\end{proof} 

The following proposition, which is an essential part of the proof of Theorem~\ref{th5.4}, is a consequence of Proposition~\ref{dos}.

\begin{proposition}\label{wuj1a}
Suppose that $\xi \in C_c^{\infty}$ is an even function such that
{\rm supp}~$\xi \subset [-1,1]$,  $\xi \ge 0$,   $\widehat{\xi}(0)=1$   and
$\widehat{\xi}^{(k)}(0)=0$  for all $1\le k  \le [s]+2$.  Then for any Borel function $ F $ with
{\rm supp} $F\subset [-N/2,N]$ we have the estimate
\begin{equation*}
	\|\delta_N(F*\xi)\|_{N,2}\le C\|\delta_NF\|_2
\end{equation*}	
and when $ s > 1/2 $ we also have
$$ \| \delta_N(F-F*\xi)\|_{N,2} \le CN^{-s}\|\delta_N F\|_{H^s}.$$

\end{proposition}

\begin{proof} 
Note that by taking $G=\delta_N F$ in \eqref{f: Np estimate  2}  we obtain the first estimate stated in this proposition. The second one concerning $ \delta_N( F -F \ast \xi) $ is proved in the Proposition \ref{dos}, see \eqref{f: EstF-Fxi}. 
\end{proof}

We start by considering the operator $ \delta_N(F-F\ast \xi)(\sqrt{\Delta_{\a,\b}}).$  We need a couple of lemmas  in order to prove Lemma~\ref{l5.4}. In what follows we let 
$ P_k^{\a,\b}f = \langle f, \psi_k^{\a,\b}\rangle_\b\, \psi_k^{\a,\b} $ stand for the projection onto the $k$-th eigenspace.

\begin{lemma}\label{chi12}
 Let $ k, N \in \mathbb{N}$ be such that $ N/2 < \lambda_k < N.$  Then for all $\a,\b\ge -1/2,$ except  for the case 
 $\beta>0$ and $\a < \b,$ we have
 \begin{equation}\label{Esrr1+2}
 \big\|\psi^{\a,\b}_k(\chi^1_N+\chi^2_{N})\big\|_{p}\le C N^{\rho(\b,p)}, 
    \end{equation}
where $ \rho(\b,p) $ is defined by
 \begin{align*}
     \rho(\b,p):=\begin{cases}
                (\frac{1}{p}-\frac{1}{2})(1+\b),\quad 1\le p\le \frac{4\b+4}{2\b+1},\\
                \frac{\b}{2}-\frac{1}{p}(1+\b),\quad \frac{4\b+4}{2\b+1}\le p\le\infty.
            \end{cases}
        \end{align*} 
Consequently,
for $\a,\b\ge -1/2$ and $1\le p\le (4\b+4)/(2\b+3)$ we have,
\begin{equation}\label{PX1X2}
\big\|P^{\a,\b}_k(\chi^1_N +\chi^2_N) \big\|_{L^p\to L^2} =\big\|\psi^{\a,\b}_k(\chi^1_N+\chi^2_{N})\big\|_{p'}\le C N^{(\b+1)(\frac{1}{p}-\frac{1}{2})-\frac12} 
\end{equation}
Moreover, when $\beta>0$ and $\a < \b$, $\big\|P^{\a,\b}_k(\chi^1_N +\chi^2_N) \big\|_{L^p\to L^2}= \infty$ for $p\ge p_*=\frac{2\b+2}{\b-\a}$ and the estimates in \eqref{PX1X2} are still valid otherwise. 
\end{lemma}
\begin{proof}
We consider the $L^{p}(d\mu_\b)$ norm of $\psi^{\a,\b}_k$ restricted to the interval $ [0, \sqrt{N/2}].$ In view of the asymptotic properties of the Laguerre functions (see \eqref{psi})
   \begin{align*}
\big\|\psi^{\a,\b}_k(\chi^1_N+\chi^2_{N})\big\|^{p}_{p}&=\int_{0}^{1/\sqrt{N}}|\lambda_k^{\a/2}x^{\a-\b}|^{p}x^{2\b+1}dx+\int_{1/\sqrt{N}}^{\sqrt{N/2}}|\lambda_k^{-1/4}x^{-(\b+\frac12)}|^{p}x^{2\b+1}dx = I_1 +I_2.
   \end{align*}
When $\a\ge\b,$ for the first integral we have
\begin{align*}
	I_1=\lambda_k^{\a p/2}\int_0^{1/\sqrt N}|x^{(\a-\b)}|^{p}x^{2\b+1}\,dx&\sim
		N^{\b(p/2-1)-1}.
\end{align*}
When $\b>0$ and $\a<\b,$ for the same integral
\begin{align*}
I_1=\lambda_k^{\a p/2}\int_0^{1/\sqrt N}|x^{(\a-\b)}|^{p}x^{2\b+1}\,dx\sim \begin{cases}
N^{\b(p/2-1)-1}, \quad  2\b+2>p(\b-\a),\\[1.2ex]
+\infty, \quad   2\b+2\le p(\b-\a).
\end{cases}
\end{align*}
\begin{align*}
  \text{And }\quad      I_2\sim \lambda_k^{-p/4} \l[\frac{x^{-p(\b+\frac12)+2\b+2}}{-p(\b+\frac12)+2\b+2}\r]_{1/\sqrt{N}}^{\sqrt{N/2}}&\le C \begin{cases}
    N^{-p/4}\l(\frac{1}{\sqrt{N}}\r)^{\b(2-p)-p/2+2}, \quad \b(2-p)\le \frac{p}{2}-2, \\
    N^{-p/4}(\sqrt{N})^{\b(2-p)-p/2+2}, \quad \b(2-p)\ge \frac{p}{2}-2
\end{cases}\\
&=C \begin{cases}
    N^{\b(p/2-1)-1}, \quad p\ge\frac{4\b+4}{2\b+1}\\
    N^{(\b+1)(1-\frac{p}{2})}, \quad p\le\frac{4\b+4}{2\b+1}.
\end{cases}
   \end{align*} 
For $\a\ge\b$, this gives the estimates in \eqref{Esrr1+2}, and for $\a<\b$, we get
\begin{equation*} \big\|\psi^{\a,\b}_k(\chi^1_N+\chi^2_{N})\big\|_{p}\le C \begin{cases}
		N^{(\frac{1}{p}-\frac{1}{2})(1+\b)},\quad 1\le p\le \frac{4\b+4}{2\b+1},\\
		N^{\frac{\b}{2}-\frac{1}{p}(1+\b)},\quad \frac{4\b+4}{2\b+1}\le p\le p^*.
	\end{cases}
\end{equation*}
Using these estimates we get
\begin{equation*} \big\|P^{\a,\b}_k(\chi^1_N +\chi^2_N) \big\|_{L^p\to L^2} =\big\|\psi^{\a,\b}_k(\chi^1_N+\chi^2_{N})\big\|_{p'}\le C N^{\b(1/2-1/p')-1/p'}=C N^{(\b+1)(1/p-1/2)-1/2}
\end{equation*}
 for $1\le p\le(4\b+4)/(2\b+3)$.
\end{proof} 

\vskip0.15truein

In what follows, it is efficient to work with $\sqrt{\Delta_{\a,\b}}$ instead of $\Delta_{\a,\b}$, because in our approach it is essential to use the corresponding wave propagator.  We can write $ S_R^\delta(\Delta_{\a,\b}) = F(\sqrt{\Delta_{\a,\b}}) $ by taking $ F(\lambda) = (1-\lambda^2/R)_+^\delta.$ Also note that $ F(\sqrt{\Delta_{\a,\b}})= G(\Delta_{\a,\b}) $ with $ G(\lambda) = (1-\lambda/R)_+^\delta.$
\begin{lemma}\label{LemChi1Chi2}
 For $N \in \mathbb{N}$ let $ F $ be a Borel function supported in $ [N/2, N].$  Then  for $\a,\b\ge -1/2$ and $1\le p\le (4\b+4)/(2\b+3)$ we have 
 \begin{equation}\label{gl}
\|F(\sqrt{\Delta_{\a,\b}})(\chi^1_{N^2}+\chi^2_{N^2})\|_{L^p \to L^2} \leq 
CN^{2(\b+1)(\frac{1}{p} - \frac{1}{2})} \|\delta_N F\|_{N^2, 2} .
\end{equation} 
\end{lemma}
\begin{proof}  For any $ f \in L^2(\R^+,d\mu_\b) $ let $ \chi_N = \chi_{N^2}^1+\chi_{N^2}^2$ and define $ g_N = \chi_N f.$ Consider
$$ \| F(\sqrt{\Delta_{\a,\b}})g_N\|_2^2= \sum_{N/2 \leq \sqrt{\lambda_k} \leq N} |F(\sqrt{\lambda_k})|^2\, \|P_k^{\a,\b} g_N\|_2^2 = \sum_{N^2/2 \leq j \leq N^2}
\sum_{ \sqrt{\lambda_k} \in [ j/N, \left(j+1)/N\right)} |F(\sqrt{\lambda_k})|^2\, \|P_k^{\a,\b} g_N\|_2^2$$
Using the estimate for the projections $ P_k \chi_N $ proved in  Lemma \ref{chi12} we obtain 
$$\sum_{ \sqrt{\lambda_k} \in [ j/N, (j+1)/N)} |F(\sqrt{\lambda_k})|^2\, \|P_k^{\a,\b} g_N\|_2^2 \leq C \sup_{\lambda \in [ j/N, (j+1)/N)} |F(\lambda)|^2  \left( \sum_{ \sqrt{\lambda_k} \in [ j/N, (j+1)/N)} 
(\sqrt{\lambda_k})^{4(\b+1)(\frac{1}{p}-\frac{1}{2})-2}\right)\,\, \|f\|_p^2.$$
In the last sum above  $ \sqrt{\lambda_k} \sim j/N \sim N,$ and the number of $ k $ such that $ j/N \leq \sqrt{\lambda_k} < (j+1)/N$ is of the order of $ j/N^2 \sim 1$ so we get the estimate
$$ \|F(\sqrt{\Delta_{\a,\b}})(\chi^1_{N^2}+\chi^2_{N^2})\|_{L^p \to L^2}^2 \leq C \, N^{4(\b+1)(1/p-1/2)} \left( \frac{1}{N^2} \sum_{N^2/2 \leq j \leq N^2} \sup_{\lambda \in [ j/N, (j+1)/N)} |F(\lambda)|^2 \right) .$$
This proves that 
$$ \|F(\sqrt{\Delta_{\a,\b}})(\chi^1_{N^2}+\chi^2_{N^2})\|_{L^p \to L^2} \leq C \, N^{2(\b+1)(1/p-1/2)} \, \| \delta_NF\|_{N^2,2}.$$
\end{proof}

The following remark is not used in the proof, as our approach relies on \eqref{3chi}. Nevertheless, it may be viewed as an interesting extension of Lemma \ref{LemChi1Chi2}.

\begin{remark}
Note that for $\a \ge -1/2$  and  $\b\ge -1/3$ we also have
\begin{equation}\label{Esrr3}
 \big\|\psi^{\a,\b}_k\chi^3_N\big\|_{p}\le C N^{\tilde{\rho}(\b,p)} \quad  \forall \, N/2< \lambda_k<2N
    \end{equation}
  where \begin{align*}
            \tilde{\rho}(\b,p):=\begin{cases}
                (\frac{1}{p}-\frac{1}{2})(1+\b),\quad 1\le p\le 4,\\
                -\frac13+\frac{1}{p}(\frac13+\b)-\frac{\b}{2},\quad 4\le p\le \frac{6\b+4}{3\b+1},\\
                \frac{\b}{2}-\frac{1}{p}(1+\b),\quad \frac{6\b+4}{3\b+1}\le p\le\infty, 
            \end{cases}
        \end{align*}  
and for $1\le p\le (6\b+4)/(3\b+3)$
\begin{equation*}
\big\|P^{\a,\b}_k\chi^3_{N} \big\|_{p\to 2} =\big\|\psi^{\a,\b}_k\chi^3_{N}\big\|_{p'}\le C N^{(\b+1)(\frac{1}{p}-\frac{1}{2})-\frac12} \quad  \forall\, N/2< \lambda_k<2N.
\end{equation*} 
When $N$ is replaced by $N^2$ here, then following a similar argument as in Lemma~\ref{chi12}, we obtain
 \begin{equation}\label{gl2}
\|F(\sqrt{\Delta_{\a,\b}})\chi^3_{N^2}\|_{L^p(d\mu_\b) \to L^2(d\mu_\b)} \le
CN^{2\tilde{\rho}(\b,p')}\|\delta_N F\|_{N^2, 2} \le 
CN^{2(\b+1)(\frac{1}{p} - \frac{1}{2})} \|\delta_N F\|_{N^2, 2} 
\end{equation} 
whenever $1\le p\le (6\b+4)/(3\b+3)$.
\end{remark}

In the following lemma we estimate the $p$ to $p$ norm of the difference operator $(F-F\ast \xi_N)(\sqrt{\Delta_{\a,\b}})$.   Note that in the first estimate the operator is self-adjoint and the term $\chi^3_{N^2}$ is included in the second estimate when $p\le 2$.

\begin{lemma}\label{l5.4}
For any Borel function $F$  supported in $ [N/2, N]$  we have
$$
\|(\chi^1_{N^2} +\chi^2_{N^2})(F-F \ast \xi )(\sqrt{\Delta_{\a,\b}})(\chi^1_{N^2} +\chi^2_{N^2})\|_{L^p \to L^p} \leq C  \|\delta_N F \|_{H^{s}}
$$
for all $s \ge \max\{1/2, 2(\b+1)|1/p-1/2|\}.$
In addition, if $p \le 2$ then 
$$
\| \chi^3_{N}(F-F \ast \xi )(\sqrt{\Delta_{\a,\b}})(\chi^1_{N^2} +\chi^2_{N^2})\|_{L^p \to L^p} \leq C  \|\delta_N F \|_{H^{s}} 
$$
for all $s \ge \max\{1/2, 2(\b+1)(1/p-1/2)\}$.
The same estimates hold for $\chi^3_{N^2}(F-F \ast \xi)(\sqrt{\Delta_{\a,\b}})\chi^1_{N^2}$.
\end{lemma}

\begin{proof} Let $ T = (F-F \ast \xi )(\sqrt{\Delta_{\a,\b}}).$ Then with $ f \in L^p(d\mu_\b), 1 \leq p \leq 2,$ applying H{\"o}lder's inequality
$$ \int_0^\infty  |(\chi^1_{N^2} +\chi^2_{N^2}) T(\chi^1_{N^2} +\chi^2_{N^2})f(x)|^p \, d\mu_\b \leq  \| T(\chi^1_{N^2} +\chi^2_{N^2})f\|_2^{2/q} \, \left( \int_0^{N/2} x^{2\b+1} dx\right)^{1/q^\prime}$$
where $ q = 2/p.$ This shows that the norm of the operator under consideration on $ L^p(d\mu_\b) $ is controlled by a constant times
$$\,N^{2(\b+1)(\frac{1}{p} - \frac{1}{2})}\|(F-F \ast \xi )(\sqrt{\Delta_{\a,\b}})(\chi^1_{N^2} +\chi^2_{N^2})\|_{L^p \to L^2}.$$ 
Under the further assumption that  \(p \le (4\beta+4)/(2\beta+3)\), we can appeal to the previous lemma and Proposition~\ref{wuj1a} to conclude the stated estimate. For other values of $ p $ we can use the self adjointness of the operator and interpolation.

\end{proof}

Having taken care of $ (F-F \ast \xi )(\sqrt{\Delta_{\a,\b}})$ it remains to treat the part $(F \ast \xi )(\sqrt{\Delta_{\a,\b}}).$ It is here we need to make use of the finite propagation speed of the wave operator $ \cos (t\sqrt{\Delta_{\a,\b}}).$ We begin with the following lemma  which is a stronger version of Lemma \ref{LemChi1Chi2}.  In what follows we let $ V(x,r) $ stand for the measure $ \mu_{\b}(B(x,r)) $ of the interval $ B(x,r) = \{y \in \R^+: |x-y| <r \}.$  We also define $ P_{B(x,r)} f= \chi_{B(x,r)} f .$ With these notations we have\\

\begin{lemma}
Assume that $ \lambda_k \sim N^2.$ Then for  $1\le p\le  \frac{4(\b+1)}{(2\b+3)}$ we have  the following estimate:
\begin{equation}\label{PB(x,r)}
\big\|P_{B(x, r)}(\chi^1_{N^2}+\chi^2_{N^2})P_k^{\a,\b} (\chi^1_{N^2}+\chi^2_{N^2})P_{B(x, r)} \big\|_{p\to p'}  \leq CV(x,
r)^{{2(\frac12}-{\frac1p})} \big( \lambda_k^{1/2} r \big)^{4(\b+1)({\frac1p}-{\frac12})} \lambda_k^{-1}
\end{equation}
for $r \geq \lambda_k^{-1/2}$, where $p'$ is the conjugate exponent of $p$.
\end{lemma}

\begin{proof} 
Consider
\begin{align*}
   \big\|P_{B(x, r)} (\chi^1_{N^2}+\chi^2_{N^2}) P^{\a,\b}_k 
(\chi^1_{N^2}+\chi^2_{N^2})P_{B(x, r)} \big\|_{p\to p'}\le \big\|\psi^{\a,\b}_k 
(\chi^1_{N^2}+\chi^2_{N^2})\chi_{B(x, r)} \big\|^2_{p'}
\end{align*}
 On the support of the function $(\chi^1_{N^2}(y)+\chi^2_{N^2}(y))\chi_{B(x, r)}(y)$ we have $ 0 \leq y \leq \sqrt{\lambda_k/2} $ and $ x-r \leq y \leq x+r.$  We first consider the case  $ 0 \leq x \leq 2r.$ We need to estimate  the integral 
$$ \left(\int_0^ {\sqrt{\lambda_k/2}} |\psi_k^{\a,\b}(y)|^{p^\prime}\, y^{2\b+1}\, dy \right)^{1/p^\prime}.$$ 
We split the integral into two parts and first estimate the integral taken over the interval $ (0, 1/\sqrt{\lambda_k}).$ Using the asymptotics of $ \psi_k^{\a,\b} $ we are led to the integral
$$ \lambda_k^{\a/2}  \left(\int_0^ {1/\sqrt{\lambda_k}}  y^{(\a-\b)p^\prime}\, y^{2\b+1}\, dy \right)^{1/p^\prime} \leq C\, \lambda_k^{\a/2}\, \lambda_k^{(-\a+\b)/2-(\b+1)(1-1/p)}$$
which is bounded by $ C\, \lambda_k^{-1/2}\, \lambda_k^{(\b+1)(1/p-1/2)}.$  By writing 
$ 1 = V(x,r)^{(1/2-1/p)}\, V(x,r)^{(1/p-1/2)} $ and using the estimate $ V(x,r) \leq C\, r^{2\b+2} $ the above can be bounded by 
$$ C \, \lambda_k^{-1/2}\, V(x,r)^{(1/2-1/p)} (r \sqrt{\lambda_k})^{(2\b+2)(1/p-1/2)} .$$ 
Next we consider the integral
$$ \left(\int_{1/\sqrt{\lambda_k}}^ {\sqrt{\lambda_k/2}} |\psi_k^{\a,\b}(y)|^{p^\prime}\, y^{2\b+1}\, dy \right)^{1/p^\prime}.$$ 
Using the asymptotics of the Laguerre functions, we are led to consider 
$$ \lambda_k^{-1/4}  \left(\int_{1/\sqrt{\lambda_k}}^ {\sqrt{\lambda_k/2}}  y^{-(\b+1/2)p^\prime}\, y^{2\b+1}\, dy \right)^{1/p^\prime}.$$ 
Under the assumption on $ p,$ we can check that $ -(\b+1/2)p^\prime+2\b+2 < 0 $ and hence the above integral is bounded by a constant multiple of 
$$ \lambda_k^{-1/4}  \, ( \sqrt{\lambda_k})^{(\b+1/2)-(2\b+2)(1-1/p)} = \lambda_k^{-1/2}\, (\sqrt{\lambda_k})^{(2\b+2)(1/p-1/2)} .$$
As before this leads to the required estimate, viz.,
$$ C \, \lambda_k^{-1/2}\, V(x,r)^{(1/2-1/p)} (r \sqrt{\lambda_k})^{(2\b+2)(1/p-1/2)}. $$ 

It remains to consider the case when $ x \geq 2r.$ In this case as $ x \geq 2r,  x-r \geq x/2 $ and $r \geq \sqrt{\lambda_k}^{-1} $ we are led to estimating 
$$ \lambda_k^{-1/4}  \left(\int_{x-r}^ {x+r}  y^{-(\b+1/2)p^\prime}\, y^{2\b+1}\, dy \right)^{1/p^\prime} \leq \lambda_k^{-1/4}  \left(\int_{x/2}^ {x+r}  y^{-(\b+1/2)p^\prime}\, y^{2\b+1}\, dy \right)^{1/p^\prime}.$$
The above  is bounded by a constant multiple of  
$$ \lambda_k^{-1/4} \, x^{-(\b+1/2)+(2\b+1)(1-1/p)}\, r^{1-1/p} \leq \lambda_k^{-1/2} \, x^{(2\b+1)(1/2-1/p)}\, r^{1/2-1/p}  ( r \sqrt{\lambda_k} )^{1/2}.$$
As $ V(x,r) $ is comparable to $ x^{2\b+1}\, r $ the above is bounded by 
$$ C \, \lambda_k^{-1/2}\, V(x,r)^{1/2-1/p} ( r \sqrt{\lambda_k} )^{1/2}.$$  Note that the condition on $ p $ implies that $ 1/2 \leq (2\b+2)(1/p-1/2) $  and hence as $ r \sqrt{\lambda_k} \geq 1,$ we get the required estimate.
\end{proof}

\begin{lemma}\label{EstProj}
Assume that $ 1 \leq p \leq \frac{4(\b+1)}{(2\b+3)}.$ Then for any Borel function $F $ with  $ {\supp} F \subset [N/2,N]$ \eqref{PB(x,r)} is equivalent to the estimate
\begin{equation}\label{spec2}
\big\|F\Big(\sqrt{\Delta_{\a,\b}}\Big) (\chi^1_{N^2}+\chi^2_{N^2})P_{B(x, r)}\big\|_{p\to 2} \leq
 CV(x,r)^{{\frac12}-{\frac1p}}(Nr)^{2(\b+1)({\frac1p}-{\frac12})}\|\delta_NF\|_{N^2,\, 2}
\end{equation}
for  all $x\in \R^+$ and $r\geq 1/N$.
For any $ \xi $ as in Proposition \ref{wuj1a}, we also have 
\begin{equation}\label{rest}
\big\|(F\ast \xi)(\sqrt{\Delta_{\a,\b}}) (\chi^1_{N^2}+\chi^2_{N^2})P_{B(x, r)}\big\|_{p\to 2} \leq CV(x,
r)^{{\frac12}-{\frac1p}} \big( Nr \big)^{2(\b+1)({\frac1p}-{\frac12})}\big\|\delta_NF\big\|_{2}. 
\end{equation} 
 \end{lemma}
\begin{proof} As in the proof of  Proposition~\ref{wuj1a} we can check that 
   $ \|\delta_N(F\ast \xi)\|_{N^2,\, 2} \le C \| \delta_NF\|_2 $ and hence  \eqref{rest} is a consequence of \eqref{spec2}. We now show that \eqref{PB(x,r)} is equivalent to \eqref{spec2}. Note that by $TT^*$ argument, 
\begin{equation*}
    \big\|P_{B(x, r)} (\chi^1_{N^2}+\chi^2_{N^2})|F|^2(\sqrt{\Delta_{\a,\b}})(\chi^1_{N^2}+\chi^2_{N^2})P_{B(x, r)} \big\|_{p\to p'}=\big\|F(\sqrt{\Delta_{\a,\b}}) (\chi^1_{N^2}+\chi^2_{N^2})P_{B(x, r)}\big\|_{p\to 2}^2.
\end{equation*}
Thus, in order to show that   \eqref{PB(x,r)} implies \eqref{spec2}, we have to check if
$$ \big\|P_{B(x, r)} (\chi^1_{N^2}+\chi^2_{N^2})|F|^2(\sqrt{\Delta_{\a,\b}})(\chi^1_{N^2}+\chi^2_{N^2})P_{B(x, r)} \big\|_{p\to p'} \leq CV(x,
r)^{{2(\frac12}-{\frac1p})} \big( Nr \big)^{4(\b+1)({\frac1p}-{\frac12})}\big\|\delta_NF\big\|^2_{N^2,2} .$$
If we let $ T_{N,r} $ stand for the operator $P_k^{\a,\b} (\chi^1_{N^2}+\chi^2_{N^2})P_{B(x, r)},$ then we can write the operator in the above inequality as
$$ \sum_{N/2 \leq \sqrt{\lambda_k} \leq N} |F(\sqrt{\lambda_k})|^2\, T_{N,r}^\ast T_{N,r}  $$
which leads to the estimate
$$ \sum_{N \leq \sqrt{\lambda_k} \leq N} |F(\sqrt{\lambda_k})|^2\, \| T_{N,r}^\ast T_{N,r}\|_{p \to p^\prime} = \sum_{N^2/2 \leq j \leq N^2}
\sum_{ \sqrt{\lambda_k} \in [ j/N, \left(j+1)/N\right)} |F(\sqrt{\lambda_k})|^2\, \| T_{N,r}^\ast T_{N,r}\|_{p \to p^\prime} .$$ 
Under the assumption that  \eqref{PB(x,r)} holds the inner sum in the above is controlled by 
$$ CV(x,r)^{{2(\frac12}-{\frac1p})} \big( N r \big)^{4(\b+1)({\frac1p}-{\frac12})} N^{-2}
\sum_{ \sqrt{\lambda_k} \in [ j/N, \left(j+1)/N\right)]} |F(\sqrt{\lambda_k})|^2. $$
Finally, as the number of $ k $ such that $ j/N \leq \sqrt{\lambda_k} < (j+1)/N$ is of the order of $ j/N^2 \sim 1$we get 
$$ \frac{1}{N^2} \sum_{N^2/2 \leq j \leq N^2}
\sum_{ \sqrt{\lambda_k} \in [ j/N, \left(j+1)/N\right)} |F(\sqrt{\lambda_k})|^2\, \leq  \frac{1}{N^2} \sum_{N^2/2 \leq j \leq N^2} \sup_{\lambda \in [j/N, (j+1)/N] }|F(\lambda)|^2 \leq \| \delta_NF\|^2_{N^2,2}.$$
For the converse we can choose $ F $ to be a bump function supported in a small neighbourhood of $ \lambda_m $ and equal to one on an even smaller neighbourhood. With such an $ F $ the inequality \eqref{spec2} gives
$$\big\|P_{B(x, r)}(\chi^1_{N^2}+\chi^2_{N^2})P_m^{\a,\b} (\chi^1_{N^2}+\chi^2_{N^2})P_{B(x, r)} \big\|_{p\to p'}  \leq C\,V(x,
r)^{{2(\frac12}-{\frac1p})} \big( N r \big)^{4(\b+1)({\frac1p}-{\frac12})}  \|\delta_NF\|_{N^2,2}^2. $$
An easy calculation shows that $  \|\delta_NF\|_{N^2,2}^2 = N^{-2} $ and hence we obtain \eqref{PB(x,r)}.
\end{proof}

\begin{remark}
    Note that Lemma \ref{EstProj} is a stronger version of Lemma \ref{LemChi1Chi2}. Indeed $$V(x,
r)^{{\frac12}-{\frac1p}} \big( Nr \big)^{2(\b+1)({\frac1p}-{\frac12})} 
=N^{2(\b+1)\left(\frac1p-\frac1s\right)}
\left(\frac{r^{2(\b+1)}}{V(x,r)}\right)^{\frac1p-\frac1s}
\le N^{2(\b+1)(\frac{1}{p} - \frac{1}{2})}. $$
\end{remark}

 Before proceeding further, let us  recall the notion of finite speed of propagation for an operator $T \colon L^p(\R^+, d\mu_\b) \to L^q(\R^+, d\mu_\b).$ 
For any $\rho > 0$ let
\begin{equation*}
 \mathcal{D}_\rho = \{ (x, y) \in \R^+ \times \R^+ : |x-y| \leq \rho \}
\end{equation*}
 stand for a neighbourhood of the diagonal in $ \R^+ \times \R^+.$ We say that $ T $ is supported in $ \mathcal{D}_\rho $ if $ \langle Tf, g \rangle = 0 $  for all $ f \in L^p, g \in L^{q^\prime} $ having supports in $ |a-x|\leq \rho_1 $ and $ |b-x|\leq \rho_2 $ with $ \rho_1+\rho_2 + \rho < |a-b|.$ If $ T $ is an integral operator with kernel $ K_T,$ then the above is equivalent to $ K_T(x,y) $ being supported in $ \mathcal{D}_\rho.$ We will be using the fact that the wave operator $ \cos (t \sqrt{\Delta_{\a,\b}}) $ is supported in $ \mathcal{D}_t$ for any  $ t>0.$\\
 
 For operators $ T $ supported in $ \mathcal{D}_\rho $ we have the following result proved in a more general context in  \cite[Lemma I.12]{COSY}.
 
 \begin{lemma}\label{I.12} Let $ T $ be a linear operator such that for all $ x \in \R^+$ and $ r >0$ the operator $ T P_{B(x,r)} $ is bounded from $ L^p(\R^+, d\mu_\b) $ to $ L^2(\R^+, d\mu_\b) $ for some $ 1 \leq p < 2.$ Assume also that $ T $ is supported in $ \mathcal{D}_\rho$ for some $ \rho >0.$ Then there exists a constant $ C= C_{p,2} $ such that
 $$ \|T\|_{p\to p} \leq C\, \sup_{x \in \R^+} \{ V(x,\rho)^{(\frac{1}{p}-\frac{1}{2})} \, \| TP_{B(x,\rho)}\|_{p\to 2} \} .$$
 \end{lemma}

We will make use of the above lemma in proving the following proposition. We closely follow the proof of \cite[Theorem I.5]{COSY} in our approach and hence we will be sketchy.

\begin{proposition}\label{PropFxi}
    Let \(N\in\mathbb N\). Suppose that \(F\) is an even Borel function  supported in $ N/2 \leq |\lambda| \leq N $ and let $ \xi $ be as in Proposition \ref{wuj1a}.
	Then there is a constant \(C\) independent of \(N\) and \(F\) such that
	\[
	\|(\chi^1_{N^2} +\chi^2_{N^2})(F \ast \xi) (\sqrt{\Delta_{\a,\b}})(\chi^1_{N^2} +\chi^2_{N^2})\|_{L^p \to L^p} \leq C  \|\delta_N F \|_{H^{s}}
	\]
	for all $s > (2\b+2)(1/p-1/2),\, 1 \leq p <2.$ 
	
\end{proposition}
\begin{proof}
	 For any even function $ G $ on $ \R ,$ Fourier inversion formula  reads as
	$$
	G(\sqrt{\Delta_{\a,\b}})
	=\frac{1}{2\pi}\int_{-\infty}^\infty \widehat{G}(t) \cos(t\sqrt{\Delta_{\a,\b}}) d t.
	$$
If we further assume that  ${\supp}\widehat{G} \subset [-r,r]$ for some $r>0 ,$ then it follows from \cite[Lemma 3]{mz} that 
	\begin{equation}\label{fsp}
		\supp  K_{G(\sqrt{\Delta_{\a,\b}})} \subseteq \mathcal{D}_r.  
	\end{equation}
	Let   $\eta \in C_c^{\infty}(\R) $   be even and  such that $\supp \eta\subseteq \{ \lambda: 1/4\leq |\lambda|\leq 1\}$ and
	$ \sum_{\ell\in \mathbb{Z}} \eta(2^{-\ell} \lambda)=1$ for all $ \lambda >0.$
	Let us set   $ \eta_\ell(\lambda) = \eta(2^{-\ell}\lambda),$ for $ \ell > 0 $ and $\eta_0(\lambda)= 1-\sum_{\ell> 0} \eta(2^{-\ell} \lambda)$  so that $ \sum_{\ell \geq 0} \eta_\ell(\lambda) =1.$  We take $ G = \delta_N F  $ and  define $ G^{(\ell)} $ as follows:
	\begin{eqnarray}\label{eq3.7}
		G^{(\ell)}(\lambda)=\frac{1}{2\pi}\int_{-\infty}^{+\infty}
		\eta_\ell(t) \hat{G}(t) \cos(t\lambda) \;dt,\,\,\,\, G_N^{(\ell)}(\lambda) = G^{(\ell)}(\lambda/N)
	\end{eqnarray}
	for $ \ell \geq 0.$ 
	In virtue of the Fourier inversion formula we have
	$$ \delta_NF(\lambda)= \sum_{\ell \ge 0} G^{(\ell)}(\lambda), \,\,\,\,  F(\lambda) =\sum_{\ell \ge 0}G_N^{(\ell)}(\lambda)$$ 
	and hence we have the decomposition
	$$ F(\sqrt{\Delta_{\a,\b}})= \sum_{\ell =0}^\infty  G_N^{(\ell)}(\sqrt{\Delta_{\a,\b}}) =\sum_{\ell =0}^\infty G^{(\ell)}(\sqrt{\Delta_{\a,\b}}/N).$$
	From the equation
	$$G^{(\ell)}(\lambda/N)= \frac{1}{2\pi}\int_{-\infty}^{+\infty}
		\eta(2^{-\ell} Nt) \hat{F}(t) \, \,\cos(t\lambda) \;dt$$
	and the finite speed of propagation of the wave operator we infer that
	$ \supp K_{G^{(\ell)}(\sqrt{\Delta_{\a,\b}}/N)} \subset \mathcal{D}_{2^{\ell}/N}. $
	
	By using the above decomposition we see that 
	$$ \| (\chi^1_{N^2} +\chi^2_{N^2}) F(\sqrt{\Delta_{\a,\b}})(\chi^1_{N^2} +\chi^2_{N^2})\|_{p\to p} $$
	$$ \leq \sum_{\ell \ge 0}\big\| G_N^{(\ell)}(\sqrt{\Delta_{\a,\b}}) (\chi^1_{N^2} +\chi^2_{N^2})\big\|_{p\to p }.$$
	In view of Lemma \ref{I.12} the above can be controlled by a constant multiple of 
	\begin{equation}\label{cos}
	\sum_{\ell \ge 0}    \sup_{x\in X}\big\{V(x,2^{\ell}/N)^{\frac1p -\frac12}\big\| G_N^{(\ell)}(\sqrt{\Delta_{\a,\b}}) (\chi^1_{N^2} +\chi^2_{N^2}) P_{B(x,2^{\ell}/N)} \big\|_{p\to 2 }\big\}.
	\end{equation}
	Since $ G_N^{(\ell)}$ are not compactly supported we cannot directly appeal to Lemma \ref{EstProj} and hence we proceed as in the proof of \cite[Theorem I.5]{COSY} using a cut-off function. We choose an even function 
	$ \zeta \in C_c^\infty(-1,1) $ supported in $ 1/8 \leq |\lambda| \leq 3/4 $ such that $ |\zeta(\lambda)| \leq 1$ and $\zeta(\lambda) = 1 $ for $ 1/4 \leq |\lambda| \leq 1/2$  and define $ \zeta_N(\lambda) = \zeta(\lambda/N).$ By  splitting $ G_N^{(\ell)} = \zeta_N\, G_N^{(\ell)} + (1-\zeta_N) G_N^{(\ell)} $ we estimate the norms of the resulting operators separately.  Note that 
	\begin{equation}\label{decomp} (\zeta_N F)\ast \xi = \sum_{\ell \geq 0} (\zeta_N G_N^{(\ell)}) \ast \xi .
	\end{equation} We will get an estimate for the part of the operator corresponding to $ (\zeta_N  F\ast \xi)(\sqrt{\Delta_{\a,\b}}).$ As in Proposition \ref{wuj1a} we will estimate the operator corresponding to $ \zeta_N F $
	
As $ \zeta_N G_N^{(\ell)} $ is compactly supported we can appeal to Lemma \ref{EstProj} to show that  norm of the operator corresponding to $ (\zeta_N G_N^{(\ell)})(\sqrt{\Delta_{\a,\b}}) $ is controlled by a constant times
	$$ V(x,2^\ell/N)^{({\frac12}-{\frac1p})} \big( N 2^{\ell}/N \big)^{2(\b+1)({\frac1p}-{\frac12})}\big\| \delta_N \zeta_N\, \delta_NG_N^{(\ell)}\big\|_{N^2, 2}.$$
	Since $ \delta_NG_N^{(\ell)} = G^{(\ell)} $ and $ |\zeta| \leq 1,$ arguing as in Proposition \ref{wuj1a} we obtain
	$$\|(\chi^1_{N^2} +\chi^2_{N^2}) (\zeta_N G_N^{(\ell)}\ast \xi)) (\sqrt{\Delta_{\a,\b}})(\chi^1_{N^2} +\chi^2_{N^2})\|_{L^p \to L^p}  $$
	$$ \leq C V(x,2^\ell/N)^{({\frac12}-{\frac1p})}  2^{2\ell(\b+1)({\frac1p}-{\frac12})}\big\|  G^{(\ell)}\big\|_{2}.$$
	Starting from \eqref{decomp}, estimating each term as above and summing up we get
	$$\|(\chi^1_{N^2} +\chi^2_{N^2}) ((\zeta_N F)\ast \xi) (\sqrt{\Delta_{\a,\b}})(\chi^1_{N^2} +\chi^2_{N^2})\|_{L^p \to L^p}  \leq C \sum_{\ell \geq 0}  2^{2\ell(\b+1)({\frac1p}-{\frac12})}\big\|  G^{(\ell)}\big\|_{ 2}.$$
Since $ G= \delta_N F, $ the sum on the right hand side of the above inequality is the norm of $ \delta_N F $ in the Besov space $ B_{2,1}^{2(\b+1)({\frac1p}-{\frac12})} $ which is bounded by $\| \delta_NF\|_{H^s}$ if 
$ s >2(\b+1)({\frac1p}-{\frac12}).$

Now we turn our attention to  estimating the operator corresponding to the part $ (1-\zeta_N)F \ast \xi.$ Since $ (1-\zeta) $ is supported in $  |\lambda| \geq 1/2,$ we  can choose $ \varphi \in C_c^\infty(1/2, 2) $ such that
$$ 1-\zeta(\lambda)  = \sum_{k=0}^\infty \varphi(2^{-k}\lambda) = \sum_{k=0}^\infty \varphi^{(k)}(\lambda) $$ 
for all $ \lambda >0.$ With  $ G^{(\ell)} $ and $ G_N^{(\ell)}$ defined as before and $ \varphi_N^{(k)}(\lambda) = \varphi^{(k)}(\lambda/N) $ we get the decomposition
$$  (1-\zeta_N)G_N^{(\ell)} = \sum_{k=0}^\infty \varphi_N^{(k)} G_N^{(\ell)}\,\,\,\,\, \sum_{\ell = 0}^\infty (1-\zeta_N)G_N^{(\ell)}\ast \xi = (1-\zeta_N)F\ast \xi.$$
In view of this decomposition we can estimate 
$$ \| ((1-\zeta_N)G_N^{(\ell)}\ast \xi)(\sqrt{\Delta_{\a,\b}})(\chi^1_{N^2} +\chi^2_{N^2}) P_{B(x,2^\ell/N)}\|_{p\to 2} $$
$$ \leq \sum_{k=0}^\infty  \| (\varphi_N^{(k)} G_N^{(\ell)}\ast \xi)(\sqrt{\Delta_{\a,\b}})  (\chi^1_{N^2} +\chi^2_{N^2}) P_{B(x,2^\ell/N)}  \|_{p\to 2}.$$ 
Since $ \varphi_N^{(k)} G_N^{(\ell)} \ast \xi$ is compactly supported each term in the above sum can be estimated to give
\begin{equation}\label{est-one}
C\, V(x,2^\ell/N)^{({\frac12}-{\frac1p})}  2^{2\ell(\b+1)({\frac1p}-{\frac12})}\big\|  \delta_N (\varphi_N^{(k)}G_N^{(\ell)}\ast \xi) \big\|_{N^2,2}.
\end{equation}
By the same argument used in Proposition \ref{wuj1a} we  have
\begin{equation}\label{est-two}
\big\|  \delta_N (\varphi_N^{(k)}G_N^{(\ell)}\ast \xi) \big\|_{N^2,2} \leq C \, \| \varphi^{(k)}  G^{(\ell)} \|_2.
\end{equation}
If $ \psi $ stands for the inverse Fourier transform of $ \eta ,$  then 
$$ \varphi^{(k)}(x)\,G^{(\ell)}(x) =  \varphi^{(k)}(x) \int_{-1}^1 \delta_NF(y) \delta_{2^{\ell}}\psi(x-y) dy.$$
Note that on the support of $ \varphi^{(k)},  |x| $ is comparable to $ 2^k.$ Therefore,  when $ x $ is in the support of $ \varphi^{(k)},$  in the integral  $ |x-y|  $ is also comparable to $ 2^k.$ Since $ \psi $ is a Schwartz class function, 
$\delta_{2^{\ell}}\psi(x-y)$ can be seen to decay like $ 2^{-M(k+\ell)} $ for any $ M \in \mathbb N.$ This gives
$$ | \varphi^{(k)}(x)\,G^{(\ell)}(x)| \leq C \varphi^{(k)}(x)\,2^{-M(k+\ell)} \|\delta_NF\|_2,\,\,\, \|\varphi^{(k)}G^{(\ell)}\|_2 \leq C 2^{k(2\b+2)}\,2^{-M(k+\ell)} \|\delta_NF\|_2.$$
The above estimate, along with \eqref{est-two} and \eqref{est-one}, leads to the estimate
$$ \| ((1-\zeta_N)G_N^{(\ell)}\ast \xi)(\sqrt{\Delta_{\a,\b}})(\chi^1_{N^2} +\chi^2_{N^2}) P_{B(x,2^\ell/N)}\|_{p\to 2}$$
$$ \leq C\, V(x,2^\ell/N)^{({\frac12}-{\frac1p})}  2^{2\ell(\b+1)({\frac1p}-{\frac12})}\, 2^{-M\ell}\, \| \delta_NF\|_2.$$
Finally, using the above in \eqref{cos} (with $ G_N^{(\ell)} $ replaced by $ (1-\zeta_N)G_N^{(\ell)} )$  we complete the proof.\\
\end{proof}

It still remains to consider the operator $\chi^3_{N^2}(F \ast \xi)(\sqrt{\Delta_{\a,\b}})\chi^1_{N^2} $ which deals with the interaction between the Bessel-type and harmonic oscillator parts.

\begin{lemma} \label{lemchi3chi1}
Under the same assumptions as in the previous proposition we have
$$ \|\chi^3_{N^2}(F \ast \xi)(\sqrt{\Delta_{\a,\b}})\, \chi^1_{N^2})\|_{L^p \to L^p}
\leq C  \|\delta_N F \|_{H^{s}}. $$
\end{lemma}
\begin{proof} With $ G = \delta_NF,$ define $ G^{(\ell)} $ and $ G_N^{(\ell)} $ as before. By defining $G_N^{loc} = \sum_{2^{\ell}/N \leq cN} G_N^{(\ell)} $ and $ G_N^\infty = \sum_{2^{\ell}/N > cN} G_N^{(\ell)} $ we can decompose 
  $ F = G_N^{loc} + G_N^\infty .$ Recall that the kernel of $ G_N^{(\ell)}(\sqrt{\Delta_{\a,\b}}) $ is supported in $ \mathcal{D}_{2^\ell/N}.$ Therefore, on the support of the operator $\chi^3_{N^2} G_N^{(\ell)}(\sqrt{\Delta_{\a,\b}})\chi^1_{N^2} $ we have $ |x-y| \leq 2^{\ell}/N, x^2> N^2/2$ and $ y^2 \leq N^2/4.$  As a consequence  $\chi^3_{N^2} G_N^{(\ell)}(\sqrt{\Delta_{\a,\b}})\chi^1_{N^2} $ is nontrivial only if $ 2^{\ell}/N > c N, c = (\sqrt{2}-1)/2$
 and hence we only need to consider $\chi^3_{N^2} G_N^{\infty}(\sqrt{\Delta_{\a,\b}})\chi^1_{N^2} .$  As $ \chi_{N^2}^3 $ is supported in $ 0 \leq x  \leq \sqrt{2}N,$ by H\"older's inequality we get 
$$ \| \chi_{N^2}^3 G_N^\infty(\sqrt{\Delta_{\a,\b}})\, \chi_{N^2}^1 \|_{p\to p} \leq C\, N^{2(\b+1)(\frac{1}{p}-\frac{1}{2})}  \|  G_N^\infty(\sqrt{\Delta_{\a,\b}})\, \chi_{N^2}^1 \|_{p\to 2}.$$
As in Proposition \ref{PropFxi}, we choose $ \zeta_N $ and consider the terms $ (\zeta_NG_N^\infty)(\sqrt{\Delta_{\a,\b}}) $ and $ ((1-\zeta_N)G_N^\infty)(\sqrt{\Delta_{\a,\b}}) $ separately. Now we can appeal to  \eqref{gl} to get
$$ \|  (\zeta_N G_N^\infty)(\sqrt{\Delta_{\a,\b}})\, \chi_{N^2}^1\|_{p\to2} \leq  C\,N^{2(\b+1)(\frac{1}{p}-\frac{1}{2})} \| \delta_N(\zeta_N G_N^\infty)\|_{N^2,2}.$$
As before this gives us the estimate
$$ \| \chi_{N^2}^3 (\zeta_NG_N^\infty \ast \xi)(\sqrt{\Delta_{\a,\b}})\, \chi_{N^2}^1 \|_{p\to p} \leq C\, N^{4(\b+1)(\frac{1}{p}-\frac{1}{2})} \,\| \delta_NG_N^\infty\|_2.$$
Recalling the definition of $ G_N^\infty $ and using the finite overlap property of $ \eta(2^\ell \cdot) $ we have
$$\| \delta_NG_N^\infty\|_2^2 \leq \sum_{2^{\ell}/N > cN} \| G^{(\ell)}\|_2^2 = \sum_{2^{\ell}/N > cN} \int_{-\infty}^\infty \eta(2^{-\ell}t)^2 |\hat{G}(t)|^2\, dt \leq C \int_{N^2}^\infty |\hat{G}(t)|^2 dt.$$
This shows that $\| \delta_NG_N^\infty\|_2 \leq C\, (1+N^4)^{-s/2} \,\| \delta_NF\|_{H^s} $ and hence for  $s>2(\b+1)(1/p-1/2)$
$$ \| \chi_{N^2}^3 (\zeta_NG_N^\infty \ast \xi)(\sqrt{\Delta_{\a,\b}})\, \chi_{N^2}^1 \|_{p\to p} \leq C\,\| \delta_NF\|_{H^s}. $$ 
With the same notations as in Proposition \ref{PropFxi} we also have
$$ \| ((1-\zeta_N)G_N^{(\ell)}\ast \xi)(\sqrt{\Delta_{\a,\b}}) \chi^1_{N^2} \|_{p\to 2}$$
$$ \leq \sum_{k=0}^\infty  \| (\varphi_N^{(k)} G_N^{(\ell)}\ast \xi)(\sqrt{\Delta_{\a,\b}})  \chi^1_{N^2}   \|_{p\to 2}.$$ 
Each term under the summation can be estimated by
$$ N^{2(\b+1)({\frac1p}-{\frac12})}\big\|  \delta_N (\varphi_N^{(k)}G_N^{(\ell)}\ast \xi) \big\|_{N^2,2} \leq C N^{2(\b+1)({\frac1p}-{\frac12})} \| \varphi^{(k)} G^{(\ell)} \|_2.$$
 Since we also have $\| \varphi^{(k)} G^{(\ell)} \|_2 \leq C 2^{k(2\b+2)}\,2^{-M(k+\ell)} \|\delta_NF\|_2$ we see that
 $$ \| ((1-\zeta_N)G_N^{\infty}\ast \xi)(\sqrt{\Delta_{\a,\b}}) \chi^1_{N^2} \|_{p\to p} \leq C N^{4(\b+1)({\frac1p}-{\frac12})}  \left( \sum_{2^{\ell}/N \geq N} \,2^{-M \ell} \right) \, \| \delta_NF\|_2.$$
 By choosing $ M $ large enough we get the required estimate.
\end{proof}

\begin{remark}\label{Rmk:chi3chi1}
Note that from Lemma~\ref{lemchi3chi1} and Lemma~\ref{l5.4}, if $p\le 2$, we obtain that for any even Borel function $F$ supported in $N/2\le |\lambda|\le N$, we have 
\begin{equation*}
	\|\chi^3_{N^2} F(\sqrt{\Delta_{\a,\b}})\chi^1_{N^2}\|_{L^p\to L^p}\le C \|\delta_N F\|_{H^s}.
\end{equation*}
If we again choose $F(\lambda)=(1-\lambda^2/N^2)^\delta_+$, then $F\in H^s$ if and only if $\delta>s-1/2$. Thus, for $S^\delta_{N^2}(\Delta_{\a,\b})=F(\sqrt{\Delta_{\a,\b}})$ and using Remark~\ref{Rmk:BRmeans}, we obtain that for $p\le 2$ and 
$$
\delta >  \max\left\{0, \quad
2(\b+1)\left|\frac{1}{2}-\frac{1}{p}\right|
-\frac{1}{2}
\right\}
$$
the  Bochner-Riesz means of order $\delta$  are uniformly bounded on $L^p
(\R^+,x^{2\b+1}dx)$, that is
$$
\sup_{N>0}\|\chi^3_{N^2}S_{N^2}^{\delta}( \Delta_{\a,\b})\chi^1_{N^2}\|_{p \to p} \le C < \infty.
$$
\end{remark}

\subsection{Conclusion}

We are now set to conclude the proof of Theorem~\ref{Th:Sufficient} that provides the sufficient condition for the uniform boundedness of Bochner-Riesz means of the operator $\Delta_{\a,\b}$. 

Note that the proof of this theorem relies on the decomposition formulated in the Splitting Lemma \ref{elem1}, by taking $T=S_N^{\delta}(\Delta_{\a,\b})$. That is, we split the operator into six parts namely, 
\begin{equation*}
	(\chi^1_N+\chi^2_N)T(\chi^1_N+\chi^2_N),\quad    (\chi^2_N+\chi^3_N)T(\chi^2_N+\chi^3_N),\quad \chi^1_N T\chi^3_N, \quad \chi^3_N T\chi^1_N,\quad  \chi^4_N T,\quad T\chi^4_N.
\end{equation*}

As mentioned in Section~\ref{sec5.1}, the last two components associated with $\chi^4_N$ are the easiest to control and are directly derived from Lemma \ref{2N}. The second component with $(\chi_N^{2}+\chi_N^3)$ is the harmonic oscillator part and its proof is a consequence of a transplantation phenomenon. It follows easily from Proposition \ref{trans} and Theorem~\ref{verlp2po}. In addition to that, for $p\le 2$, the term $\chi^1_N T\chi^3_N$ can be absorbed into the component $(\chi^2_N+\chi^3_N)T(\chi^2_N+\chi^3_N)$, and the boundedness of this part follows from Corollary~\ref{TT}.

We studied the first component with $(\chi_N^{1}+\chi_N^{2})$ in Section~\ref{sec5.1}. This is a Bessel-type part of our operator as it corresponds to Bessel operators. The approach to handling this operator is based on the finite speed of propagation property for the corresponding wave equation. Hence, we needed to replace our operator $\Delta_{\a,\b}$ by $\sqrt{\Delta_{\a,\b}}$ and $N$ in the characteristic functions by $N^2$. Here we followed the spectral multiplier approach and considered the operator $F(\sqrt{\Delta_{\a,\b}})$, which eventually gives the uniform boundedness of $S^\delta_N(\Delta_{\a,\b})$ by taking $F(\lambda)= (1-\lambda^2/N)^\delta_+$.

We split this operator $F(\sqrt{\Delta_{\a,\b}})$ into $(F*\xi)\big(\sqrt{\Delta_{\a,\b}}\big)+(F-F*\xi)\big(\sqrt{\Delta_{\a,\b}}\big)$
where $\xi$ is a smooth, compactly supported function on $\RR$ satisfying the assumptions stated in the Proposition~\ref{wuj1a}. Here, the second part $(F-F*\xi)\big(\sqrt{\Delta_{\a,\b}}\big)$ is easier to handle and the $L^p$ to $L^p$ estimates for this part are obtained in Lemma~\ref{l5.4}. The first part is where the finite propagation speed of the corresponding wave operator is used and Proposition \ref{PropFxi} provides the desired estimates.

Note that for the $\chi^3_N T\chi^1_N$ component as well, our approach is based on the same finite speed of propagation property. Therefore, by the same splitting of the operator  $F(\sqrt{\Delta_{\a,\b}})$ as before, we prove the result for this component in Lemma~\ref{l5.4} for the $F-F*\xi$ part and in Lemma~\ref{lemchi3chi1} for the $F*\xi$ part. The corresponding results for the uniform boundedness of the Bochner-Riesz means $S^\delta_N(\Delta_{\a,\b})$ for these two components are mentioned in Remark~\ref{Rmk:chi1chi2} and Remark~\ref{Rmk:chi3chi1}.

This concludes the proof of Theorem~\ref{Th:Sufficient}.\\

\noindent
\begin{center}{\bf Acknowledgements}
\end{center}

H. Sharma thanks Dorothee Frey for being her host during her two-year postdoctoral research at KIT, where this project was initiated. She acknowledges the financial support from the Alexander von Humboldt Foundation during her stay at KIT and for her research visit to Macquarie University to work on this project. A. Sikora also thanks Dorothee Frey for hosting and supporting his two research visits to Karlsruhe to work on this project. A. Sikora and H. Sharma were also supported by Australian Research Council  Discovery Grant  DP200101065. The work of S. Thangavelu is supported by the ANRF research grant NSC/2025/000054. He wishes to thank ANRF, IISc and UQ for providing him with comfortable chairs to ponder over the problem studied in this work.

\medskip

%

\bibliographystyle{abbrv}
\bibliography{Bibliography}
\end{document}